\documentclass[12pt]{article}

\usepackage{amsmath,amsthm}
\usepackage{amssymb}
\usepackage{epsfig}

\usepackage{easybmat}
\usepackage{graphicx,hyperref}
\usepackage{url}

\usepackage{cite}

\newtheorem{theorem}{Theorem}[section]
\newtheorem{lemma}[theorem]{Lemma}

\newtheorem{cor}[theorem]{Corollary}

\theoremstyle{definition}
\newtheorem{define}[theorem]{Definition}
\newtheorem{remark}[theorem]{Remark}

\newcommand{\dsty}{\displaystyle}
\newcommand{\casefrac}[2]{ {\textstyle{ \frac{#1}{#2} } } } 

\newcommand{\pf}{\noindent {\bf Proof: }}
\newcommand{\enpf}{\hfill $\Box$ \vspace{.2in} }

\begin{document}

\title{Morse Theory and Relative Equilibria \\ in the Planar $n$-Vortex Problem}

\author{Gareth E. Roberts\thanks{
Dept. of Mathematics and Computer Science,
College of the Holy Cross,
groberts@holycross.edu}
}

\maketitle

\begin{abstract}
Morse theoretical ideas are applied to the study of relative equilibria in the planar $n$-vortex problem.
For the case of positive circulations, we prove that the Morse index of a critical point of the Hamiltonian 
restricted to a level surface of the angular impulse is equal to the number of pairs of real eigenvalues of the corresponding relative equilibrium periodic solution. 
The Morse inequalities are then used to prove the instability of some families of relative equilibria in the four-vortex problem with
two pairs of equal vorticities.   We also show that, for positive circulations, relative equilibria cannot accumulate on the collision set.
\end{abstract}


\vspace{.1in}


\section{Introduction}

In the study of the planar $n$-vortex problem, periodic solutions for which the configuration of vortices rotates rigidly about the
center of vorticity play a crucial role.  Such solutions are known as {\em relative equilibria}, since they are fixed points in a rotating
coordinate system.   Rigidly rotating vortex configurations, sometimes referred to as {\em vortex crystals}~\cite{aref-newton}, 
are frequently observed in physical experiments (e.g., rotating superfluid $^4$He~\cite{yarmchuk} or
Bose-Einstein condensates~\cite{navarro}) as well as in numerical models of natural phenomena (e.g., 
the eyewalls of hurricanes~\cite{davis, KossSchub}).  
Consequently, investigating the stability of relative equilibria is of great significance, 
not only for identifying stable solutions, but also for understanding the local structure of the flow
in a neighborhood of the periodic solution.

The purpose of this paper is to use Morse theoretic ideas to study the stability of relative equilibria in the planar
$n$-vortex problem.   Let $H$ denote the Hamiltonian for the problem and $I$ the angular impulse.  
Relative equilibria are found as critical points of the smooth function $H$ restricted to
the level surface $I = I_0$ (an ellipsoid).  In~\cite{g:stability}, it was shown that for same-signed circulations, a relative equilibrium $z$ is linearly stable if and only if
$z$ is a nondegenerate minimum of $H$ restricted to $I = I_0$.   As $I$ is a conserved quantity for the point vortex problem, 
linear stability actually implies non-linear stability as well, due to a theorem of Dirichlet's.  We extend these results by establishing a connection
between the Morse index of a relative equilibrium and its eigenvalues in a rotating frame.  Specifically, we show that the Morse index is equivalent
to the number of pairs of real eigenvalues $\pm \lambda$.   This implies that the index of a relative equilibrium is directly related to
its degree of instability.   Our result corresponds well with the recent work of Barutello, Jadanza, and Portaluri, who studied the instability of relative
equilibria for mechanical systems with homogeneous and logarithmic potential functions~\cite{BJP}.
The authors show that an odd Morse index implies the relative equilibrium is linearly unstable.\footnote{Note
that the main result in~\cite{BJP} does not directly apply to our problem because the mechanical systems treated therein include
kinetic energy (phase space of dimension~$4n$), a quantity absent from the point vortex problem.}

We also demonstrate the effectiveness of the Morse inequalities by applying them to a special case of the four-vortex problem for which
the number and type of relative equilibria are known.   The case considered here is the {\em two equal pairs problem}, where the
circulations are $\Gamma_1 = \Gamma_2 = 1$ and $\Gamma_3 = \Gamma_4 = m$, with $m$ a parameter in the interval $(0,1)$.   
In~\cite{HRS}, the number of relative equilibria for
each $m \in (0,1)$ was rigorously shown to be~34.  Moreover, the configurations were grouped into three classes: 
convex (6), concave (16), and collinear (12).
A configuration is called {\em concave} if one vortex is located strictly inside the convex hull of the other three;
{\em convex} if no vortex is contained in the convex hull of the other three, and
{\em collinear} if all vortices lie on a common line.   The convex solutions for the two equal pairs problem are minima~\cite{g:stability}, 
while the collinear configurations are known to be unstable with index~2.  Using the Morse inequalities, we quickly deduce that
the 16 concave solutions have index~1 and are thus unstable.  The simplicity of the calculation underscores the power of the Morse
theoretical approach, although establishing the nondegeneracy of the concave solutions is a challenging computation that
requires explicit formulas for the configurations in terms of~$m$.

Our Morse theoretical approach builds on the work of Smale~\cite{smaleProbs}, Palmore~\cite{palmoreCW, palmoreClass, palmoreColl, palmore}, 
and particularly Moeckel, whose recent book chapter on central configurations
in the $n$-body problem gives a clear and comprehensive treatment concerning the topological approach to
the study of relative equilibria~\cite{rick-book}.

 In the next section, we define relative equilibria for the planar $n$-vortex problem and
provide the topological framework to handle the inherent symmetries of the problem.  We explain how to calculate the Morse index using
a special modified Hessian matrix.  We also prove the analog of Shub's lemma from celestial mechanics, showing that, for a 
fixed choice of positive circulations, relative equilibria are bounded away from the collision set.   In Section~3 we review the relevant
theory from~\cite{g:stability} concerning the linear stability of relative equilibria and establish the explicit connection 
between the index and the number of real eigenvalues.  
Section~4 focuses on applying the Morse inequalities to the two equal pairs problem.  In this section we use techniques
from computational algebraic geometry to derive formulas for the concave solutions and prove that $H$ is a Morse function.  
Symbolic and Gr\"{o}bner bases computations were performed with Magma~\cite{magma} and Maple\textsuperscript{TM}~\cite{maple}.
Eigenvalues were computed using Matlab~\cite{matlab} to check our results.  Figures \ref{fig:KiteSetup} and~\ref{fig:Asymm} 
were created with the open-source software Sage~\cite{sage}.

\section{Relative Equilibria}
\label{sec:rel-eq}

The equations of motion for $n$ planar point vortices form a first-order Hamiltonian system, as outlined
by Kirchhoff~\cite{kirchhoff}.  We let $z_i \in \mathbb{R}^2$ denote the position of the $i$th vortex and 
$r_{ij} = \|z_i - z_j\|$ represent the distance between the $i$th and $j$th vortices.
Let $z = (z_1, \ldots, z_n) \in \mathbb{R}^{2n}$ represent the vector of positions.
Each vortex has a {\em circulation} or {\em vorticity} denoted by $\Gamma_i$, where $\Gamma_i$ is a 
nonzero real number.  The Hamiltonian function for the point vortex problem is 
$$
H(z) = -\sum_{i<j} \Gamma_i \Gamma_j \ln (r_{ij}) .
$$
The motion of the $i$th vortex is then given by
\begin{equation}
\Gamma_i  \dot z_i \;  = \;  J \frac{\partial H}{\partial z_i}  \; = \;   J  \sum_{j \neq i}^n \frac{\Gamma_i \Gamma_j}{r_{ij}^2}(z_j - z_i), \quad 1 \leq  i \leq n ,
\label{eq:vorts}
\end{equation}
where 
$
J=
\begin{bmatrix}
0 & 1\\
-1 & 0
\end{bmatrix}
$
is the standard $2 \times 2$ symplectic matrix.  \\[0.04in]

Define $\Gamma = \sum_i  \Gamma_i $ as the {\em total circulation} of the system.    We will assume throughout that $\Gamma \neq 0$.
The {\em center of vorticity}, given by $c =  \frac{1} {\Gamma} \sum_i \Gamma_i z_i$,  is thus well-defined.  This is the analogue of the
center of mass in the $n$-body problem.

\subsection{Relative equilibria as critical points}

A relative equilibrium is a periodic solution of~(\ref{eq:vorts}) where the configuration of vortices rigidly rotates about~$c$.
\begin{define}
A {\em relative equilibrium} is a solution of~(\ref{eq:vorts}) of the form 
\begin{equation}
z_i(t) = c + e^{- \omega J t} (z_i(0) - c) , \quad \mbox{ for each } i \in \{1, \ldots , n\},
\label{eq:DefRE}
\end{equation}
that is, a uniform rotation with angular velocity~$\omega \neq 0$ about the center of vorticity~$c$.
\end{define}

Locating relative equilibria is a challenging algebra problem.  The initial positions $z = (z_1, \ldots , z_n)$ of a relative equilibrium must satisfy the following system
of equations:
\begin{equation}
\sum_{j \neq i}^n \frac{\Gamma_i \Gamma_j}{r_{ij}^2} (z_j - z_i)  + \omega \, \Gamma_i (z_i  - c)  \; = \;  0,  \quad  i \in \{1, \ldots, n\}.
\label{eq:rel-equ}
\end{equation}
If the position vector $z \in \mathbb{R}^{2n}$ is a solution of equation~(\ref{eq:rel-equ}), then it is standard practice to refer to $z$ as a relative
equilibrium, with the understanding that $z$ generates a periodic solution in the form of equation~(\ref{eq:DefRE}).  
We follow this practice here.  In celestial mechanics, the distinction is made clearer 
by referring to $z$ as a {\em central configuration}.

There are several symmetries present in the $n$-vortex problem that must be accounted for.  
If $z$ is a solution of equation~(\ref{eq:rel-equ}), it is straight-forward to check that the following transformations of $z$
(scaling, translation, and rotation, respectively) are also relative equilibria:
\begin{itemize}
\item[{\bf (i)}]  $\kappa z = (\kappa z_1, \ldots, \kappa z_n)$ for any scalar $\kappa > 0$  ($c \mapsto \kappa c, \omega \mapsto \omega/\kappa^2$),

\item[{\bf (ii)}]  $z - \zeta = (z_1 - \zeta, \ldots, z_n - \zeta)$ for any $\zeta \in \mathbb{R}^2$ ($c \mapsto c - \zeta$, same $\omega$),

\item[{\bf (iii)}]  $Az = (Az_1, \ldots, Az_n)$, where $A \in \mbox{SO}(2)$ ($c \mapsto Ac$, same $\omega$).
\end{itemize}
Thus, relative equilibria are not isolated.   We deal with these symmetries by specifying a scaling and center of vorticity~$c$, 
and identifying solutions that are equivalent under a rotation.  When counting or classifying solutions, 
we view relative equilibria as members of an equivalence class.  Note that reflecting each vortex about a coordinate
axis also produces a relative equilibrium (with a new~$c$ and the same~$\omega$), but 
these will be treated as distinct solutions.


The quantity 
$$
I(z) \;  = \;  \sum_{i=1}^n \, \Gamma_i \|z_i - c\|^2 ,
$$
called the {\em angular impulse with respect to the center of vorticity}, 
measures the size of the system with respect to~$c$.
It is the analogue of the {\em moment of inertia} in the $n$-body problem.   Using equation~(\ref{eq:vorts}), it is straight-forward to check that
$I$ is an integral of motion for the planar $n$-vortex problem~\cite{newton}.

One advantage of the angular impulse is that it facilitates a topological approach to the study of relative equilibria.
Note that system~(\ref{eq:rel-equ}) can be written more compactly as
\begin{equation}
\nabla H(z) + \frac{\omega}{2} \nabla I(z) \; = \; 0 ,
\label{eq:cc}
\end{equation}
where $\nabla$ is the usual gradient operator.
Thus, relative equilibria (regarded as points in~$\mathbb{R}^{2n}$) are critical
points of the Hamiltonian $H$ restricted to a level surface of $I$, where
the constant $\omega/2$ can be regarded as a Lagrange multiplier.

An important identity involving the Hamiltonian is
\begin{equation}
\nabla H(z) \cdot z \; = \; -L ,  \quad  \mbox{where } L \; = \;  \sum_{i < j} \Gamma_i \Gamma_j  
\label{eq:homogH}
\end{equation}
and $\cdot$ represents the standard Euclidean inner product.
Taking the dot product on both sides of equation~(\ref{eq:cc}) with $z - c = (z_1 - c, \ldots, z_n - c)$ and applying~(\ref{eq:homogH}) yields
the useful formula
\begin{equation}
\omega \; = \;  \frac{L}{I(z)} .
\label{eq:omega}
\end{equation}
The constant $L$, known as the {\em total vortex angular momentum}, plays an important role
in the study of relative equilibria in the planar $n$-vortex problem.   For same-signed circulations, we always have $L > 0$,
while for mixed signs, it is possible to obtain $L \leq 0$.  When $L$ changes sign, bifurcations in stability typically occur.
For the special case $L = 0$, we necessarily have $I = 0$ since $\omega \neq 0$ is assumed.  
Any relative equilibrium with $L = 0$ is always degenerate~\cite{g:stability}.

\subsection{A Morse theoretical approach}

To apply Morse theory to the study of relative equilibria, we follow Moeckel's approach in the $n$-body setting,
as expounded in~\cite{rick-book}.   Unless otherwise stated, we will always assume that the circulations satisfy $\Gamma_i > 0 \; \forall i$.
In order to eliminate the translational invariance and fix the scaling, we
restrict to the {\em normalized configuration space}
$$
{\cal N} \; = \;   \{ z \in \mathbb{R}^{2n}: c = 0, I(z) = 1 \} .
$$
The manifold ${ \cal N}$ is diffeomorphic to the sphere ${\bf S}^{2n - 3}$.
If $z$ is a critical point of $H|_{\cal N}$, then it can be shown that $z$ satisfies equation~(\ref{eq:cc}) and is thus a relative equilibrium.

Let $\Delta = \{z  \in \mathbb{R}^{2n}: z_i = z_j \, \mbox{ for some } i \neq j\}$ be the {\em collision set} containing all configurations where
two or more vortices coincide.   Since the
function $H(z)$ is undefined on $\Delta$, we will work away from this set.
To eliminate the rotational symmetry, we define the quotient space
$$
{\cal M} \; = \;  ({\cal N} - \Delta)/\mbox{SO}(2)
$$
of dimension $2n - 4$.   Since SO(2) $\cong {\bf S}^1$ acts freely on ${\cal N} - \Delta$, the quotient space ${\cal M}$ is also a manifold.
A relative equilibrium will be called {\em nondegenerate} if it is a nondegenerate critical point of $H$ restricted to~${\cal M}$.  
A criterion in terms of eigenvalues will be given in Definition~\ref{Def:RelEquDeg}.
Rigorously verifying that a particular relative equilibrium is nondegenerate can be a difficult task, 
as demonstrated with the four-vortex examples discussed in Sections \ref{SubSec:Kites} and~\ref{SubSec:Asymm}.   
Assuming that all of its critical points are nondegenerate, 
we can regard $H$ as a Morse function on ${\cal M}$.

For the case of three vortices, ${\cal N}/\mbox{SO}(2)$ is diffeomorphic to ${\bf S}^2$ and is known as the {\em shape sphere}, since
it represents the space of all triangles up to translation, scaling, and rotation.  The manifold~${\cal M}$ is thus the shape sphere minus the three points corresponding
to binary collisions (triple collision is eliminated because $I(z) = 1$).  
It is well known that $H|_{\cal M}$ has five critical points, all nondegenerate:
two equilateral triangles at the North and South Poles (minima) and three collinear configurations on the equator (saddles).
The shape sphere was a useful framework for proving the existence of the famous
figure-eight orbit in the three-body problem~\cite{cm}.

\subsubsection{The Hessian and modified Hessian}

Suppose that $z \in {\cal N}$ is a relative equilibrium for a fixed choice of circulations.  
Since $I(z) = 1$, we see that $\omega = L$ by formula~(\ref{eq:omega}).  
This means that {\em all} relative equilibria in the normalized configuration space will rotate with the {\em same} frequency and in the same direction.  
The {\em Morse index} of $z$ is the dimension of the largest subspace of $T_z{\cal N}$ for which the Hessian quadratic form is negative definite.
This is equivalent to the number of negative eigenvalues of a matrix representation of the Hessian.  
It is easier to compute the Morse index using matrices of size $2n \times 2n$, rather
than working in local coordinates on ${\cal M}$.

Introduce the function $G(z) = H(z) + (L/2) I(z)$ and let $D^2H(z)$ denote the symmetric matrix of second partial derivatives of~$H$.   
A critical point of $H|_{\cal N}$ is also a critical point of $G$.  Since
$H|_{\cal N}$ and $G|_{\cal N}$ differ by a constant, their Hessians are identical on $T_z{\cal N}$.  This justifies
the following definition.

\begin{define}
The {\em Hessian} of $H|_{\cal N}$ at a relative equilibrium $z$ is the restriction of the $2n \times 2n$ matrix
$$
D^2G(z) \; = \;  D^2H(z) + L M
$$
to the tangent space of ${ \cal N}$ at $z$, where
$M = {\rm diag}\{ \Gamma_1, \Gamma_1, \ldots, \Gamma_n, \Gamma_n \} $. 
\end{define}

Assuming that $\Gamma_i > 0 \; \forall i$, the matrix $M$ is positive definite 
and we can define the ``mass'' inner product 
$$
< \! v,w \! > \; = \;  v^T M w .
$$
In order to compute the Morse index of~$z$, it is easier to work with the {\em modified Hessian}
$$
M^{-1} D^2G(z) \; = \;  M^{-1} D^2H(z) + L I ,
$$
where $I$ now represents the $2n \times 2n$ identity matrix.   The matrix $M^{-1} D^2G$ is symmetric with respect to $< \! \ast, \ast \! >$ and therefore has only
real eigenvalues.  
Using Sylvester's Law of Inertia~\cite{CarlMeyer}, the number of negative (or zero) eigenvalues of the Hessian and the modified
Hessian are the same.

\begin{lemma}
Let $z \in {\cal N}$ be a relative equilibrium with $\Gamma_i > 0 \; \forall i$.  Then the Morse index of $z$ is equal to the number of negative eigenvalues of the
modified Hessian $M^{-1} D^2H(z) + L I .$
\end{lemma}

\pf
Let $\delta_i = [[\frac{i+1}{2}]]$, where $[[\ast]]$ represents the greatest integer function.  Define the vectors $b_i = (1/\sqrt{\Gamma_{\delta_i}}) e_i$, $i \in \{1, \ldots, 2n\}$ 
where $\{e_i\}$ are the standard basis vectors of $\mathbb{R}^{2n}$.  The vectors $\{b_i\}$ form an orthonormal basis of $\mathbb{R}^{2n}$ with respect
to $< \! \ast, \ast \! >$.   Let $P$ be the $2n \times 2n$ matrix whose $i$th column is $b_i$.  Then $P$ satisfies $P^T M P = I$ or $P^{-1} = P^T M$.

Recall that two square matrices $A_1$ and $A_2$ are called congruent if there exists an invertible matrix $Q$ such that $A_2 = Q^T A_1 Q$.
Since
$$
P^{-1} (M^{-1}D^2H(z) + LI) P \; = \;  P^TD^2H(z)P + LI   \; = \;  P^T(D^2H(z) + LM) P,
$$
we see that $M^{-1}D^2H(z) + LI$ is similar to $P^T(D^2H(z) + LM)P$, which in turn, is congruent to $D^2H(z) + LM$.  
Sylvester's Inertia Law states that the number of negative eigenvalues is identical for congruent matrices.  Thus,
the modified Hessian $M^{-1} D^2G(z)$ and the Hessian $D^2G(z)$ have the same number of negative eigenvalues.
\enpf

\subsubsection{Trivial Eigenvalues}
\label{subsec:trivialEvals}

Define the $2n \times 2n$ block-diagonal matrix $K = {\rm diag}\{ J, J, \ldots, J \} $, 
where $J=\begin{bmatrix}0&1\\-1&0\end{bmatrix}$.  
Let $z_i = (x_i, y_i)$ 
and define the quantities
$a_{ij} = (y_i - y_j)^2 - (x_i - x_j)^2$ and $b_{ij} = -2(x_i - x_j)(y_i - y_j)$.
If we write
$$
D^2H(z) \; = \;  
\begin{bmatrix}
A_{11} & A_{12} & \cdots & A_{1n} \\
\vdots &    &    &   \vdots \\
A_{n1} &  A_{n2}  & \cdots &  A_{nn}
\end{bmatrix},
$$
it is straight-forward to check that the off-diagonal blocks ($i \neq j$) are given by
$$
A_{ij}  \;  = \;  
\frac{\Gamma_i \Gamma_j}{r_{ij}^4} 
\begin{bmatrix}
a_{ij} & b_{ij} \\[0.1in]
b_{ij} & -a_{ij}
\end{bmatrix} ,
$$
while the diagonal blocks satisfy 
$
A_{ii} \; = \;  - \sum_{j \neq i} A_{ij}.
$
One crucial property of $D^2H(z)$ is that it anti-commutes with the matrix $K$:
\begin{equation}
D^2H(z) K \; = \;  -K D^2H(z) .
\label{eq:anti-commute}
\end{equation}
This follows quickly by observing that $J$ anti-commutes with each $A_{ij}$.

The modified Hessian $M^{-1} D^2H(z) + LI$ always has the eigenvalues $L, L, 2L,$ and $0$ corresponding to the symmetries discussed earlier.
We will refer to these eigenvalues as {\em trivial}.  
Due to conservation of the center of vorticity, the vectors $s = [1, 0, 1, 0, \ldots, 1, 0]^T$ and $Ks$ are in the kernel of
$M^{-1} D^2H(z)$. This follows directly from the block diagonal structure of $D^2H(z)$. 
Thus, the modified Hessian has the eigenvalue $L > 0$ repeated twice.

The other trivial eigenvalues arise from the scaling and rotational symmetries.  To see this, differentiate identity~(\ref{eq:homogH})
with respect to~$z$.  Then, assume that $z$ is a relative equilibrium and substitute in equation~(\ref{eq:cc}).  This yields the relation
\begin{equation}
M^{-1} D^2H(z) \, z \; = \;  L z .
\label{eq:eigRE}
\end{equation}
Hence the relative equilibrium itself, regarded as a vector in $\mathbb{R}^{2n}$, is an eigenvector of the modified Hessian
with eigenvalue $2L$.   Using property~(\ref{eq:anti-commute}), the vector $Kz$ is in the kernel of both the Hessian and
the modified Hessian.  This reflects the fact that relative equilibria are not isolated on~${\cal N}$.

The vectors $s, Ks$, and $z$ are all orthogonal to $T_z({\cal N})$ with respect to our inner product, and are thus ignored when
computing the Morse index.  The zero eigenvalue arising from the rotational symmetry is accounted for when reducing to the
quotient manifold~${\cal M}$.

\begin{define}
For any relative equilibrium~$z \in {\cal N}$, the modified Hessian always has the four trivial eigenvalues $L, L, 2L,$ and $0$, where $L = \sum_{i < j} \Gamma_i \Gamma_j$.
If the remaining $2n - 4$ eigenvalues are nonzero, then $z$ is a {\em nondegenerate} critical point of $H$ restricted to ${\cal M}$.  
\label{Def:RelEquDeg}
\end{define}

\subsection{Shub's lemma in the vortex setting}

Recall that the quotient manifold ${\cal M}$ excludes the collision set $\Delta$ and is therefore non-compact.  
To handle this issue, we now show that, for a fixed choice of positive circulations, the critical points
of $H|_{\cal N}$ are bounded away from~$\Delta$.
The analogous result in the $n$-body problem is called Shub's lemma~\cite{shub}.

Moeckel gives a nice argument for Shub's lemma in~\cite{rick-book} (see Prop.~2.8.7).  
Surprisingly, neither Shub's original proof nor Moeckel's extends to the $n$-vortex setting,
essentially because of identity~(\ref{eq:homogH}).   Shub's proof breaks down since the derivative he computes no longer approaches infinity;
Moeckel's argument fails to generalize because the angular velocity $\omega$ is always constant.
Nevertheless, our proof incorporates ideas from both arguments.

\begin{theorem}
For a fixed choice of circulations $\Gamma_i > 0$,
there is a neighborhood of~$\Delta$ in ${\cal N}$ which contains no relative equilibria.
\label{thm:ShubVorts}
\end{theorem}

\pf
Suppose the claim was false.  Then, since ${\cal N}$ is compact, there would exist a sequence of relative equilibria $\{z^k \} \subset {\cal N}$ converging
to some configuration $z' \in {\cal N} \cap \Delta$ on the diagonal.  By equation~(\ref{eq:cc}), we have
\begin{equation}
\nabla H(z^k) \; = \;  -L M z^k
\label{eq:ShubPf}
\end{equation}
for each~$k$.  
We will compute the directional derivative of $H$ for a well-chosen direction to show that the magnitude of the left-hand side of 
equation~(\ref{eq:ShubPf}) is unbounded as $k \rightarrow \infty$.  Since the right-hand side of equation~(\ref{eq:ShubPf}) is clearly bounded, 
we have a contradiction.

Without loss of generality, we may group the vortices so that that the first $l$ vortices ($l \geq 2$) are all approaching the same point $z_1'$, while the remaining vortices
are bounded away from~$z_1'$.  In other words, $\dsty{\lim_{k \rightarrow \infty} z_i^k = z_1'}$ for each $1 \leq i \leq l$ and
$\dsty{\lim_{k \rightarrow \infty} z_i^k = z_i' \neq z_1'}$ for each $l+1 \leq i \leq n$.
Consider the vector $v^k = (v_1^k, v_2^k, \dots, v_n^k)$ defined as
$$
v_i^k = \left\{ \begin{array}{cc}
                 z_1' - z_i^k   &   \mbox{for $1 \leq i \leq l$,} \\[0.07in]
                 0  & \mbox{for $l+1 \leq i \leq n$} 
              \end{array}
        \right.
$$
and let $u^k = v^k/||v^k||$ be the unit vector in the direction of $v^k$.  Note that $||v^k|| \rightarrow 0$ as $k \rightarrow \infty$.
The vector $u^k$ corresponds to a perturbation that directs all vortices in the first cluster toward their common limiting point.  

We compute the directional derivative of $H$ in the direction of $u^k$ by grouping all pairs of indices $(i, j)$ with $i < j \leq l$ together.
We find that
\begin{eqnarray}
\nabla H(z^k) \cdot u^k  & = &   \sum_{i=1}^l  \frac{\partial H}{\partial z_i} (z^k) \cdot \frac{v^k_i}{||v^k||} \nonumber  \\
& = &  \frac{1}{||v^k||}  \left[   \sum_{i < j}^l  \frac{\Gamma_i \Gamma_j}{r_{ij}^2} ( (z_j^k - z_i^k) \cdot (z_1' - z_i^k) + (z_i^k - z_j^k) \cdot (z_1' - z_j^k)) \right. \nonumber \\[0.07in]
&    &  \qquad   \left.  +  \;  \sum_{i=1}^l  \sum_{j = l+1}^n  \frac{\Gamma_i \Gamma_j}{r_{ij}^2} (z_j^k - z_i^k) \cdot (z_1' - z_i^k) \vphantom{\frac12}\right]  \label{eq:doubleSum}  \\
& = &   \frac{1}{||v^k||}  \left[    \sum_{i < j}^l \Gamma_i  \Gamma_j  + {\cal F} (z^k)  \right]  \nonumber
\end{eqnarray}
where ${\cal F} (z^k)$ is the double sum on line~(\ref{eq:doubleSum}).  Since the vortices $l+1 \leq i \leq n$ are bounded away from the first cluster in the limit, we see that
${\cal F} (z^k)$ approaches 0 as $k \rightarrow \infty$.  It follows that $\dsty{\lim_{k \rightarrow \infty} \nabla H(z^k) \cdot u^k = \infty}$, which implies
that $||\nabla H(z^k)||$ is unbounded as $k \rightarrow \infty$.
\enpf

When the circulations have opposite signs, ${\cal N}$ is no longer compact and 
it is possible for relative equilibria to accumulate on the collision set $\Delta$.  For example, 
if $\Gamma_1 = \Gamma_2 = \Gamma_3 = \Gamma_4 = 1$ and $\Gamma_5 = -1/2$, there exists a continua of
relative equilibria where the four equal vortices are positioned at the vertices of
a rhombus with the fifth vortex located at its center~\cite{g:continuum}.   The angle between any two adjacent sides can
be used to parametrize the continuum.  As this angle approaches zero, two vortices on a diagonal of the rhombus approach the central vortex,
limiting on a triple collision.  Despite such anomolous counterexamples, if the vortex angular momentum for a given cluster
of vortices does not vanish, then it is not possible for that cluster to limit on collision.

\begin{theorem}
Consider the $n$-vortex problem with mixed-sign circulations and assume that $L \neq 0$.  
Suppose that $\{z^k \} \subset {\cal N}$ is a sequence of relative equilibria converging to some point $z' \in \Delta$
and let ${\cal S}$ be a subset of the indices $\{1, \ldots, n\}$ corresponding to a cluster of vortices approaching collision.  
Then
$$
L'  \; = \;   \sum_{\substack{ i < j \\[0.02in] i, j \in {\cal S} } }  \Gamma_i  \Gamma_j \; = \;  0.
$$
\label{thm:CollVorts}
\end{theorem}

\pf
Since $L \neq 0$, equation~(\ref{eq:ShubPf}) still holds.  In this case, we do not want to obtain a contradiction as $k \rightarrow \infty$,
so $|| \nabla H(z^k) ||$ must remain bounded in the limit.  But if $L' \neq 0$, then the computation of the directional derivative in the proof of Theorem~\ref{thm:ShubVorts} 
would imply that $|| \nabla H(z^k) ||$ becomes unbounded.  Thus, $L' = 0$ is required for any subsets of vortices colliding in the limit.  
\enpf

\begin{remark}
\begin{enumerate}
\item  Consider the $1+$rhombus continuum discussed above and suppose that vortices 1 and~2 are on one diagonal while vortices 3 and~4 are on the other. 
Then $L' = \Gamma_1 \Gamma_2 + \Gamma_1 \Gamma_5 + \Gamma_2 \Gamma_5 = 0$ and 
$L' = \Gamma_3 \Gamma_4 + \Gamma_3 \Gamma_5 + \Gamma_4 \Gamma_5 = 0$, in accordance with the theorem.

\item  Theorem~\ref{thm:CollVorts} also applies to families of relative equilibria with {\em changing} vorticities.    If a family of relative equilibria has some
subset of vortices approaching collision, then the limiting values of the $\Gamma_i$ in this subset must satisfy $L' = 0$ (or else $L = 0$ in the limit).  
This condition is easy to guarantee if the $\Gamma_i$ approach zero;
however, examples exist for non-vanishing circulations as well.  For example, there exists a family of four-vortex collinear relative
equilibria with circulations $\Gamma_1 = \Gamma_2 = \Gamma_3 = 1$ and $\Gamma_4 = m$ that limit on triple collision as $m \rightarrow -1/2$~\cite{gBrian}.
In general, values of parameters for which $L'$ vanish are likely candidates for bifurcations.  

\end{enumerate}
\end{remark}

\section{Linear Stability and the Morse Index}

We now focus on the connection between the Morse index of a relative equilibrium and the linear stability
of the corresponding periodic solution, assuming throughout that $\Gamma_i > 0 \; \forall i$.  
We first review the salient points on linear stability
from~\cite{g:stability}.   The key idea is to exploit the fact that $D^2H(z)$ and $K$ anti-commute.

\subsection{Linear stability of relative equilibria}

Suppose that $z \in {\cal N}$ is a relative equilibrium.  By definition, this means that
the center of vorticity $c$ is at the origin, the angular impulse $I(z)$ equals unity and the angular
velocity $\omega$ is equal to the constant~$L$.  The simplest way to approach the dynamical stability 
of~$z$ is to use rotating coordinates and treat $z$ as a rest point of the corresponding flow.

Recall that $M = {\rm diag}\{ \Gamma_1, \Gamma_1, \ldots, \Gamma_n, \Gamma_n \} $ and
that $K$ is the $2n \times 2n$ block-diagonal matrix
containing $J$ on the diagonal.   From~\cite{g:stability}, the matrix that determines the linear stability of $z$ is given by
$$
B \; = \;  B(z) \; = \;  K( M^{-1} D^2 H(z) + L I ) ,
$$
where $I$ is the $2n \times 2n$ identity matrix.  
We will refer to $B$ as the {\em stability matrix}.  
It is the linearization of the planar $n$-vortex problem in rotating coordinates about~$z$.
Since the system is Hamiltonian, the characteristic polynomial of $B$ is even and the eigenvalues come in pairs $\pm \lambda$.  
For~$z$ to be linearly stable, the eigenvalues must lie on the imaginary axis.

Let $p(\lambda)$ denote the characteristic polynomial of the stability matrix $B$.
Equation~(\ref{eq:anti-commute}) yields a factorization of $p(\lambda)$ in terms of the eigenvalues
of $M^{-1} D^2 H(z)$.   If $v$ is an eigenvector of $M^{-1} D^2 H(z)$ with eigenvalue $\mu$, 
then $Kv$ is also an eigenvector with eigenvalue $- \mu$.  Consequently, 
$\mbox{span}\{ v, K v\}$ is an invariant subspace of~$B$, which yields
$\lambda^2 + L^2 - \mu^2$ as a factor of $p(\lambda)$.   Since all circulations are assumed to
be positive, the matrix $M^{-1} D^2H(z)$ is symmetric with respect to $< \! \ast, \ast \! >$ and thus has a full set of orthonormal eigenvectors
$\{v_j, Kv_j\}$ with corresponding eigenvalues $\pm \mu_j$.  
It follows that $p(\lambda)$ factors completely into $n$ even, quadratic polynomials:
\begin{equation}
p(\lambda) \; = \;  \lambda^2 (\lambda^2 + L^2) \prod_{j = 1}^{n-2} (\lambda^2 + L^2 - \mu_j^2).
\label{eq:charPoly}
\end{equation}

The first two factors in~(\ref{eq:charPoly}) arise from the symmetries discussed in Section~\ref{subsec:trivialEvals}.
From identity~(\ref{eq:eigRE}), we see that $z$ is an eigenvector of $M^{-1} D^2H(z)$ with eigenvalue $\mu = L$,
yielding the factor $\lambda^2$ and repeated zero eigenvalues (illustrating the fact that relative
equilibria are not isolated rest points).  Similarly, the vector $s = [1, 0, 1, 0, \ldots, 1, 0]^T$ 
is in the kernel of $M^{-1} D^2H(z)$, producing
the factor $\lambda^2 + L^2$ and the eigenvalues $\pm i L$.  These eigenvalues reflect the conservation of
the center of vorticity (invariance under translation).

Linear stability is defined by working over the appropriate subspace of $\mathbb{R}^{2n}$.
Set $V = \mbox{span} \{z, Kz \}$ and let $V^\perp$ denote the orthogonal complement of $V$ with respect to $< \! \ast, \ast \! >$, that is,
$$
V^\perp \; = \;  \{ w \in  \mathbb{R}^{2n} : w^T M v = 0 \; \; \forall v \in V \} .
$$
The invariant subspace $V$ accounts for the two zero eigenvalues.  
The vector space $V^\perp$ has dimension $2n - 2$ and is invariant under $B$.  
We also have that $V \cap V^\perp = \{0\}$ since $L > 0$ (see Lemma~2.6 in~\cite{g:stability}).  
This leads to the following definition for linear stability.

\begin{define}
For a relative equilibrium $z \in {\cal N}$, the stability matrix $B$ always has the four trivial eigenvalues $0, 0, \pm i L$.  
We call $z$ {\em nondegenerate} if the remaining $2n - 4$ eigenvalues are nonzero.  
A nondegenerate relative equilibrium  is {\em spectrally stable} if the nontrivial 
eigenvalues lie on the imaginary axis, and {\em linearly stable} if, in addition, the restriction of $B$ to $V^\perp$ 
has a block-diagonal Jordan form with blocks
$
\begin{bmatrix}
0 &  \beta \\  
-\beta & 0 \\
\end{bmatrix}.
$
\end{define}

\begin{remark}
The meaning of nondegeneracy in this context is consistent with that of Definition~\ref{Def:RelEquDeg}
since a vector is in the kernel of the modified Hessian if and only if it is also in the kernel of $B$.
Working over $V^\perp$ to determine linear stability is analogous to restricting to $T_z({\cal N})$ to calculate the index of $z$.
\end{remark}

From equation~(\ref{eq:charPoly}), we see that $z$ is linearly stable if and only if $|\mu_j| < L$ for each 
nontrivial eigenvalue $\mu_j$ of $M^{-1} D^2H(z)$ (the trivial eigenvalues being $0, 0, \pm L$).   Due to the special
factorization of the characteristic polynomial, spectral and linear stability are actually equivalent concepts.
The only way to lose stability is for the eigenvalues to become zero (an additional degeneracy) and then form a real pair $\pm \lambda$.
The situation is more complicated for circulations of mixed sign.  In this case, $M$ is no longer positive definite and the matrix
$M^{-1} D^2H(z)$ may have complex eigenvalues, leading to quartic factors of the characteristic polynomial and eigenvalues
of the form $\pm \alpha \pm i \beta$ (see Lemma~2.5, part~(b) in~\cite{g:stability}).

\subsection{Relating the Morse index to the eigenvalues of $B$}

Next we identify the specific connection between the index of a relative equilibrium and the eigenvalues of the corresponding periodic
solution.   Our notation and key matrices are summarized in Table~\ref{table:notation} for the reader's convenience.  

\renewcommand{\arraystretch}{1.65}
\begin{table}[h]
\begin{center}
\begin{tabular}{|c|c|c|}
\hline 
{\bf Name}  &     \makebox[2.8in]{{\bf Matrix}}      &    {\bf Eigenvalue Symbol} \\
\hline
        		               &  $M^{-1}D^2H(z)$                           		       &  $\mu$      \\
Stability Matrix       &   $B(z) =  K( M^{-1} D^2 H(z) + L I)$                   &   $\lambda$   \\
Hessian                   &  $D^2G(z) =  D^2H(z) + L M$    		       &       \\
Modified Hessian  &  $M^{-1} D^2G(z) = M^{-1}D^2H(z) + L I$         &    $\nu$    \\
\hline 
\end{tabular} 
\end{center}
\caption{The key matrices for calculating the index and eigenvalues of a relative equilibrium
$z \in {\cal N}$.
$L$ is the total vortex angular momentum (a positive constant) and $I$ is the $2n \times 2n$ identity matrix.}
\label{table:notation}
\end{table}
\renewcommand{\arraystretch}{1.0}

The Morse index of $z$ is equal to the number of negative eigenvalues of $M^{-1} D^2G(z)$ while the linear stability of~$z$
is determined by the eigenvalues of $M^{-1} D^2H(z)$.  Since the difference of these two matrices
is a scalar multiple of the identity matrix, it is straight-forward to compare their
eigenvalues.

\begin{lemma}
Suppose that $\nu$ is an eigenvalue of $M^{-1} D^2G(z)$ with eigenvector $v$.  Then,

\begin{itemize}
\item[{\bf (i)}]  $2L - \nu$ is an eigenvalue of $M^{-1} D^2G(z)$ with eigenvector $Kv$, and

\item[{\bf (ii)}]   $\mu = \nu - L$ is an eigenvalue of $M^{-1} D^2H(z)$.
\end{itemize}
\label{Lemma:eigConnect}
\end{lemma}

\pf
{\bf (i)} Using identity~(\ref{eq:anti-commute}), we have the following sequence of implications:
\begin{eqnarray*}
(M^{-1} D^2H(z) + L I) v \; = \;  \nu  v  & \Longrightarrow  &  (K M^{-1} D^2H(z) + L K) v \; = \;  \nu K v \\
&  \Longrightarrow &  (- M^{-1} D^2H(z) + L I) Kv \;  = \;  \nu K v \\
& \Longrightarrow &  (M^{-1} D^2H(z) - L I) Kv \; = \;   -\nu K v \\
& \Longrightarrow &  (M^{-1} D^2H(z) + L I) Kv \; = \;   (2L - \nu) K v \\
& \Longrightarrow &  M^{-1} D^2G(z) K v  \; = \;   (2L - \nu) K v .
\end{eqnarray*}

{\bf (ii)}  Since $(M^{-1} D^2H(z) + L I) v  =   \nu v$, we have $M^{-1} D^2H(z) v = (\nu - L) v$.
\enpf

Item {\bf (i)} of Lemma~\ref{Lemma:eigConnect} shows that the eigenvalues of the modified Hessian 
$M^{-1} D^2G(z)$ come in pairs of the form
$(\nu_j, 2L - \nu_j)$.  This yields a simple upper bound for the Morse index of $z$, denoted as~ind$(z)$.
This bound was first given by Palmore~\cite{palmore}.  The same upper bound is attained in the planar $n$-body problem
(\cite{palmoreClass, rick-book}).

\begin{theorem}
Suppose that $\Gamma_j > 0 \; \forall j$ and that $z$ is a relative equilibrium.   Then
$$
{\rm ind}(z) \; \leq \; n - 2,
$$
with equality holding whenever $z$ is a collinear configuration.
\end{theorem}

\pf
Let $\nu_j$ be an eigenvalue of $M^{-1} D^2G(z)$.
If $\nu_j < 0$, then $2L - \nu_j > 0$.  Hence, at most half of the eigenvalues of $M^{-1} D^2G(z)$
can be negative and thus
$$
{\rm ind}(z) \; \leq  \;  \frac{2n - 4}{2} \; = \;  n - 2 .
$$

The fact that any collinear $n$-vortex relative equilibrium has index $n-2$ is stated by Palmore in~\cite{palmore}, although
no proof is given.  It can be verified by generalizing a cunning topological argument due to Conley from the $n$-body problem
(see \cite{pacella} or \cite{rick-book} for an explanation of this argument).
\enpf

Due to equation~(\ref{eq:charPoly}), the nontrivial eigenvalues of the stability matrix $B$ 
are of the form $\pm \sqrt{\mu_j^2 - L^2}$, where $\mu_j$ is a nontrivial eigenvalue 
of~$M^{-1} D^2H(z)$.   Since $\mu_j$ is always real, the eigenvalues solely consist of
a real pair $\pm \lambda_j$, a pair of zero eigenvalues (degenerate case), or a pure imaginary pair $\pm i \beta_j$.
If $\nu_j$ is an eigenvalue of the modified Hessian
$M^{-1} D^2G(z)$, then item {\bf (ii)} of Lemma~\ref{Lemma:eigConnect} implies that 
\begin{equation}
\pm \sqrt{(\nu_j - L)^2 - L^2} \; = \;  \pm \sqrt{ \nu_j(\nu_j - 2L)}
\label{eigStabMat}
\end{equation}
are eigenvalues of the stability matrix $B$.  This leads to one of our main results.

\begin{theorem}
Suppose that $\Gamma_i > 0 \; \forall i$ and that $z$ is a relative equilibrium.
The Morse index of $z$ is equal to the number of pairs of
real (nonzero) eigenvalues $\pm \lambda_j$ of the corresponding periodic solution.
\label{Thm:Main}
\end{theorem}

\pf
Since translating or scaling $z$ does not change the index nor the eigenvalue structure,
we can assume that $z \in {\cal N}$.  
The nontrivial eigenvalues of the modified Hessian come in pairs of the form $(\nu_j, 2L - \nu_j)$, $j \in \{1, \ldots, n-2\}$. 
Since the values in each pair are equidistant from $L > 0$, we
may assume, without loss of generality, that $\nu_j \leq L$ for each~$j$.   

There are three possible outcomes for the eigenvalues
of the stability matrix based on the sign ($+, -,$ or $0$) of $\nu_j$.  First, if $\nu_j < 0$, then the quantity under the radical in formula~(\ref{eigStabMat})
is positive, and we obtain a pair of real eigenvalues of $B$ of the form $\pm \lambda_j$.  
Secondly, if $\nu_j = 0$, then the relative equilibrium is degenerate with a pair of zero eigenvalues.  Finally, if $0 < \nu_j \leq L$, then
the quantity under the radical in formula~(\ref{eigStabMat}) is negative, and we obtain a pure imaginary pair 
of eigenvalues of $B$ of the form $\pm i \beta_j$.   Thus, the only way to obtain a real (nonzero) pair of eigenvalues for the stability matrix is to have a negative eigenvalue
of the modified Hessian.  It follows that the Morse index of~$z$, which is equivalent to the number of negative eigenvalues of the modified Hessian, is
{\em precisely} the number of real (nonzero) pairs $\pm \lambda_j$ of eigenvalues of the stability matrix.  
\enpf

\begin{cor}
Suppose that $\Gamma_j > 0 \; \forall j$ and that $z$ is a relative equilibrium.
Then $z$ is linearly stable if and only if it is a nondegenerate minimum of $H$ subject
to the constraint $I = I_0$.
\end{cor}

\pf
This fact is the main result in~\cite{g:stability}.   If $z$ is linearly stable, then it cannot have any real pairs of eigenvalues.  By Theorem~\ref{Thm:Main},
the Morse index is zero and $z$ must be a (nondegenerate) minimum.  Conversely, if $z$ is a nondegenerate minimum, then $0 < \nu_j \leq L$
for each $j \in \{1, \ldots, n-2\}$, where $\nu_j$ is a nontrivial eigenvalue of the modified Hessian.   Spectral stability (and therefore linear stability)
now follows from formula~(\ref{eigStabMat}).
\enpf

\begin{remark}
\begin{enumerate}

\item  Theorem~\ref{Thm:Main} reveals a direct relationship between the Morse index and the instability of the relative equilibrium since
the index is equivalent to the number of real, positive eigenvalues.  There is also a connection between ``unstable'' vectors in $V^\perp$ and
directions in $T_z({\cal N})$ that decrease $H|_{\cal N}$.  To see this, suppose that $v$ is an eigenvector of $M^{-1}D^2H(z)$ with a negative eigenvalue $\mu$ satisfying
$\mu < -L$.   Then $\{v, Kv\}$ is an invariant subspace of $B$ that yields a pair of real eigenvalues $\pm \lambda$.  But we also have that
$$
v^T D^2G(z) v \;  = \;  (\mu + L) ||v||^2  \;  < \;   0,
$$ 
where the norm is computed with respect to the mass inner product.
This shows that the value of $H|_{\cal N}$ decreases in the direction of the eigenvector~$v$.

\item  Theorem~\ref{Thm:Main} is valid even if $z$ is degenerate.   Every negative eigenvalue of the modified Hessian corresponds to a real pair of
eigenvalues of the stability matrix (and vice-versa), regardless of the number of zero eigenvalues.

\item  Note that $\nu_j(2L - \nu_j) \leq L^2$ for any value of $\nu_j$.  
If $\pm i \beta_j$ are eigenvalues of the stability matrix, then formula~(\ref{eigStabMat}) shows that $|\beta_j| \leq L = \omega$.  
In other words, the angular velocity for each component in the center manifold of the linearized flow is always less than or equal to the angular
velocity of the relative equilibrium itself.  This is true whether $z$ is stable or not.

\end{enumerate}
\label{remark:mainThm}
\end{remark}

\section{The Morse Inequalities}

In this section we apply the Morse inequalities and Theorem~\ref{Thm:Main} to determine the linear stability of two families of relative equilibria
in the four-vortex problem.  The advantage of this approach is that it gives a quick argument for the instability of the solutions, avoiding
the need to compute the eigenvectors and eigenvalues of $M^{-1} D^2H(z)$ directly.

Recall that a relative equilibrium $z$ is a critical point of the Hamiltonian $H$ restricted to the manifold 
${\cal M} \; = \;  ({\cal N} - \Delta)/\mbox{SO}(2)$.   For a fixed choice of positive circulations, Theorem~\ref{thm:ShubVorts}
enables us to work on a compact space away from the singular set~$\Delta$.
The Morse inequalities relate the indices of the
critical points to the topology of~${\cal M}$ and can be written in polynomial form as
$$
\sum_k  \gamma_k t^k  \; = \;  \sum_k  p_k t^k + (1 + t) Q(t),
$$
where $\gamma_k$ is the number of critical points of index~$k$, $p_k$ is the $k$th Betti number of~${\cal M}$
(the rank of the homology group $\tilde{H}_k({\cal M}, \mathbb{R})$), and
$Q(t)$ is a polynomial with non-negative integer coefficients.

The polynomial $P(t) = \sum_k  p_k t^k$ is called the {\em Poincar\'{e} polynomial}.  Since we have removed the
collision set~$\Delta$ from our space, the topology of~${\cal M}$ is nontrivial.  
Using induction, Moeckel derives the following formula for the Poincar\'{e} polynomial of ${\cal M}$ for the planar 
$n$-body problem~\cite{rick-book}.  Since the topology of the planar $n$-vortex problem is identical, this polynomial
is valid in our setting as well.

\begin{theorem}
For the planar $n$-vortex problem, the Poincar\'{e} polynomial for ${\cal M} \; = \;  ({\cal N} - \Delta)/\mbox{SO}(2)$ is
\begin{equation}
P(t) \; = \;  (1 + 2t)(1+ 3t) \cdots (1 + (n-1)t) .
\label{eq:Poincare}
\end{equation}
\end{theorem}

Recall that for $n=3$, the manifold ${\cal M}$ is the shape sphere minus three points, which is diffeomorphic to the plane
with two points removed.  Thus, the Betti numbers are $p_0 = 1$ and $p_1 = 2$ yielding $P(t) = 1 + 2t$.  This concurs with
formula~(\ref{eq:Poincare}) when $n=3$.

Consider the planar four-vortex problem with circulations $\Gamma_1 = \Gamma_2 = 1$ and $\Gamma_3 = \Gamma_4 = m$, where
$m \in (-1, 1]$ is a parameter.  In~\cite{HRS}, an exact count on the number and type of relative equilibria solutions 
is determined in terms of~$m$ (see Table~1 in~\cite{HRS}).  For the case $0 < m < 1$, there are $34$
distinct relative equilibria:  6 convex configurations (isosceles trapezoid, rhombus),
16 concave configurations (kites, asymmetric), and 12 collinear configurations.  
In~\cite{g:stability}, it is shown that the convex configurations are each linearly stable.  
Here we will prove that the concave configurations have a Morse index of~1 and are thus linearly unstable.
Animations of these configurations for varying~$m$
can be found at \url{http://mathcs.holycross.edu/~groberts/Research/vort-movies.html}.

In order to apply the Morse inequalities, we first need to confirm that $H|_{\cal M}$ is a Morse function.  This requires checking that
the concave relative equilibria for $0 < m < 1$ are nondegenerate.  In order to accomplish this, we first find analytic formulas for each
kind of solution.

\subsection{Two pairs of equal-strength vortices: kites}
\label{SubSec:Kites}

In this section we derive formulas for the kite families of relative equilibria where vortices three and four lie on the axis of symmetry.  
We use rectangular coordinates and Gr\"{o}bner bases with symmetric coordinates to find analytic expressions for the positions as
a function of the parameter~$m$.   For an excellent source on Gr\"{o}bner bases and invariant group theory, see~\cite{CLO}.

\begin{figure}[tb]
\begin{center}
\includegraphics[width=263bp]{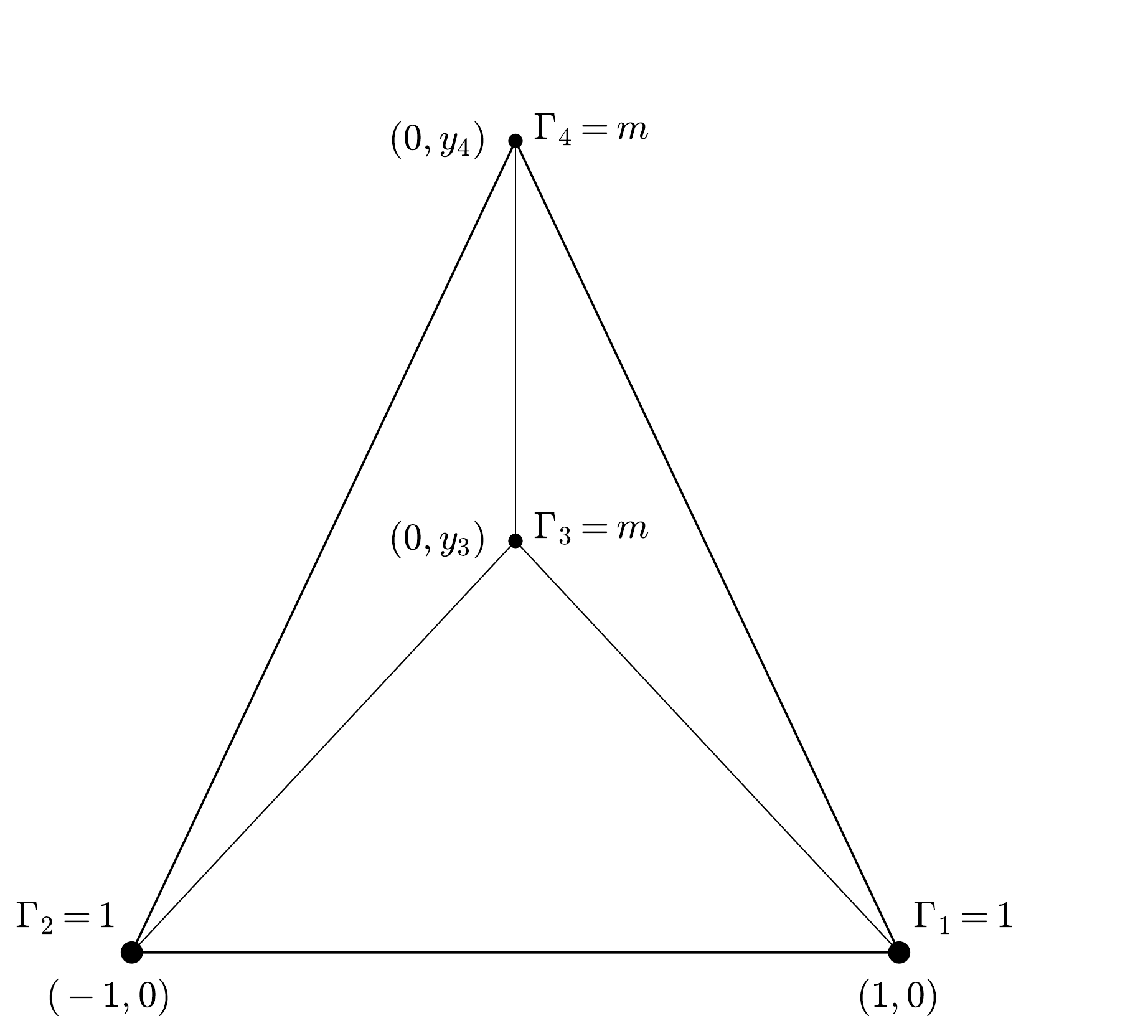}
\hspace{-0.25in}
\includegraphics[width=263bp]{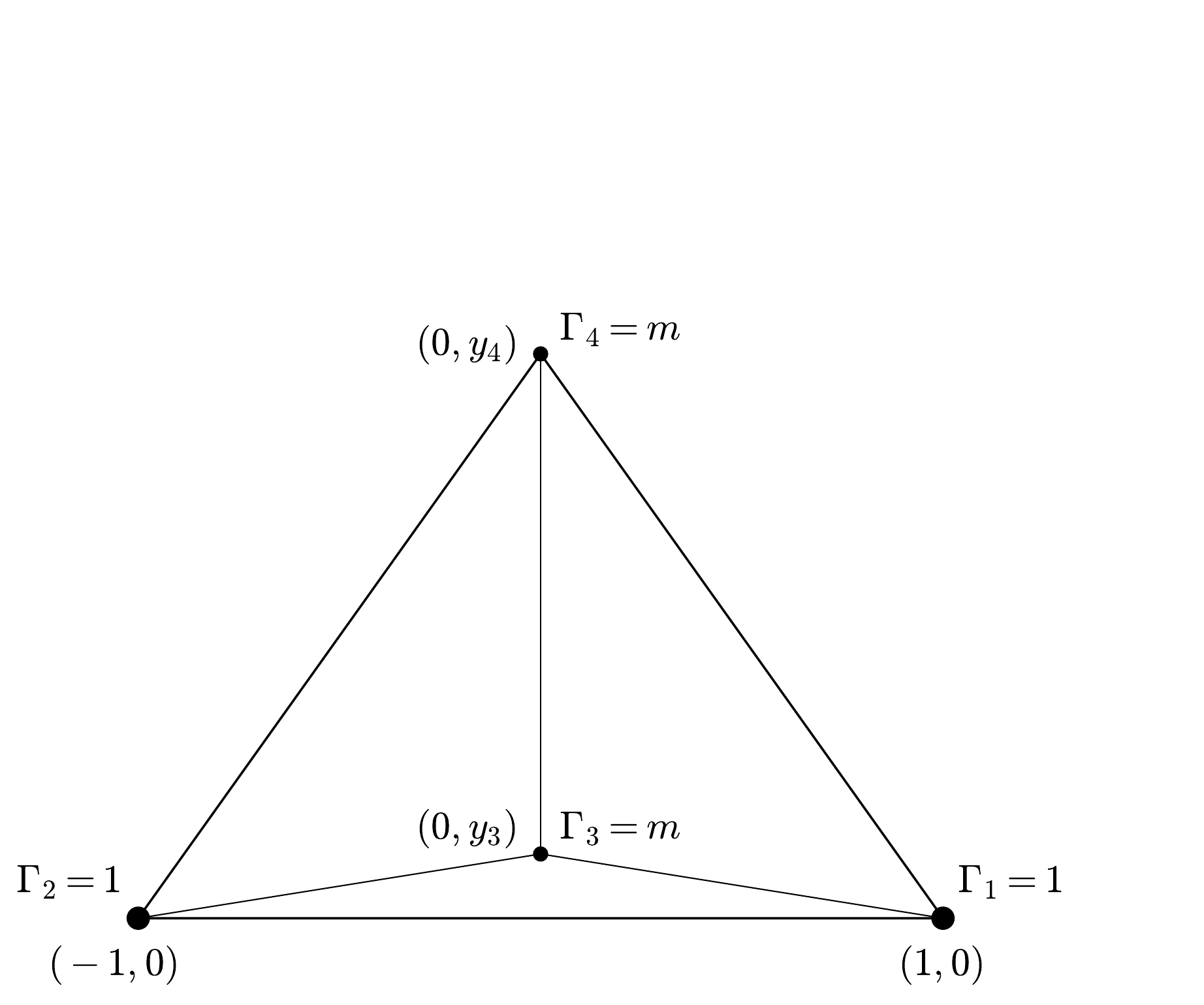}
\end{center}
\vspace{-0.2in}
\caption{Two geometrically distinct concave kite solutions for the parameter value $m = 0.6$.}
\label{fig:KiteSetup}
\end{figure}

We position the vortices at $z_1=(1,0), z_2 = (-1,0), z_3=(0,y_3)$, and $z_4=(0,y_4)$, where $y_3$ and $y_4$ are unknown
(see Figure~\ref{fig:KiteSetup}).
If $y_3 = -y_4$, the configuration forms a rhombus, a case that has already been studied in great detail in Section~4.2 of~\cite{g:stability}
and Section~7.4 of~\cite{HRS}.   The rhombus family that exists for $0 < m < 1$ is linearly stable.  We will assume that $y_3 + y_4 \neq 0$.
The center of vorticity is $c = (0, (m(a+b))/(2m+2))$.  Although $z \not \in {\cal N}$, it can easily be rescaled and translated into that space,
transformations that do not effect the nondegeneracy of~$z$.

Substituting our coordinates into system~(\ref{eq:rel-equ}) yields three independent equations:
\begin{eqnarray}
-\frac{1}{2} - \frac{m}{y_3^2 + 1} - \frac{m}{y_4^2 + 1} + \omega  & = &  0,  \label{eq:kite1} \\[0.07in]
\frac{y_3}{y_3^2+1} + \frac{y_4}{y_4^2+1} -  \frac{\omega(y_3 + y_4)}{2m + 2}  & = &  0,   \label{eq:kite2}  \\[0.07in]
-\frac{2y_3}{y_3^2+1}  - \frac{m}{y_3 - y_4} + \frac{\omega(m (y_3 - y_4) + 2y_3)}{2m + 2} & = &  0.  \label{eq:kite3}
\end{eqnarray}
Each of equations (\ref{eq:kite1}) through~(\ref{eq:kite3}) can easily be solved for~$\omega$.  Equating the first pair of expressions for $\omega$ and
the last pair produces the following polynomial system:
\begin{eqnarray}
(y_3^2 + 1)(y_4^2 + 1) - 4(y_3 y_4 + 1) + 2m(y_3 - y_4)^2 &  = &  0, \label{eq:kite4}  \\
2y_3 y_4 (y_3 - y_4)^2 + m(y_3^3 y_4 - 3y_3^2 y_4^2 + y_3 y_4^3 - 2y_3 y_4 - 1) & = &  0.  \label{eq:kite5}
\end{eqnarray}
Any solution $(y_3, y_4, m)$ satisfying both equations (\ref{eq:kite4}) and~(\ref{eq:kite5}) will yield a kite relative equilibrium
with angular velocity~$\omega$ given by~(\ref{eq:kite1}).

Note that equations (\ref{eq:kite4}) and~(\ref{eq:kite5}) possess the symmetry $(y_3, y_4) \mapsto (y_4, y_3)$, a consequence 
of $\Gamma_3 = \Gamma_4$.   To simplify the computations, we introduce the symmetric variables $\sigma$ and $\rho$ defined by
\begin{eqnarray}
\sigma & = &  y_3 + y_4,  \label{eq:sigma}  \\
\rho  & = &  y_3 y_4  , \label{eq:rho}  
\end{eqnarray}
and use Gr\"{o}bner bases to eliminate $y_3$ and $y_4$.  Specifically, we
compute a Gr\"{o}bner basis for equations (\ref{eq:kite4}) through~(\ref{eq:rho}) with respect to the lex order
$y_4 > y_3 > \sigma > \rho > m$.  We also saturate with respect to $y_3 - y_4$ to eliminate solutions where
vortices three and four collide.  Two of the polynomials in the resulting basis are
\begin{eqnarray}
\rho^2(m+2) + \rho(2m^2 - 2m - 6) + 2m^2 + m   & = &  0,  \label{eq:rhoM}\\
\sigma^2(2m + 1) + \rho^2 - (8m + 6)\rho - 3 & = &  0.  \label{eq:sigmaM}
\end{eqnarray}
Since equations (\ref{eq:rhoM}) and~(\ref{eq:sigmaM}) are quadratic in $\rho$ and $\sigma$, the problem can now be completely solved
in terms of the parameter~$m$.

\begin{remark}
We note the distinct advantage of using symmetric coordinates;  computing a Gr\"{o}bner basis
for equations (\ref{eq:kite4}) and~(\ref{eq:kite5}) with respect to the lex order $y_4 > y_3 > m$ produces a complicated eighth-degree
polynomial in $y_3$ with coefficients in~$m$.  Moreover, the sign of the variable~$\rho$ determines the type of configuration,
with $\rho > 0$ yielding a concave kite configuration and $\rho < 0$ corresponding to a convex kite configuration.  
\end{remark}

\begin{theorem}
Let $\rho$ and $\sigma$ be defined by
\begin{eqnarray}
\rho & = &  \frac{-m^2 + m + 3 \pm  \sqrt{(m^2 - 1)(m^2 - 4m -9)} }{m+2} , \label{eq:rhoForm} \\[0.07in]
\sigma^2 & = &  \frac{2[ (5m^2 + 10m+3)\rho  +m^2 + 2m + 3]}{(2m+1)(m+2)}.  \nonumber
\end{eqnarray}
If $y_3$ and $y_4$ are chosen to be distinct roots of the
quadratic $y^2 - \sigma y + \rho$, then $z = (1,0, -1,0, 0, y_3, 0, y_4)$ gives a kite relative equilibrium
with vorticities $\Gamma_1 = \Gamma_2 = 1$ and $\Gamma_3 = \Gamma_4 = m$.   
The number and type of these solutions are given as follows:

\begin{itemize}
\item[{\bf (i)}]  At $m = 1$, there are four concave solutions consisting of an equilateral triangle with a vortex at the center.
These solutions are degenerate with nullity equal to three.

\item[{\bf (ii)}]  For $0 < m < 1$, there are eight concave solutions corresponding to 
two geometrically distinct kite families.  As $m$ passes through 0, one of these families persists smoothly and becomes convex, 
continuing for $-1/2 < m < 0$.  At $m = 0$, this family has three collinear vortices and the exterior triangle 
becomes equilateral.

\item[{\bf (iii)}]  Let $\widetilde{m} \approx -1.6804$ represent the only real root of the cubic $5m^3 + 7m^2 + 3m + 9$.
For $-2 < m < \widetilde{m}$, there are four convex solutions corresponding to one geometrically distinct family of 
kites.  These kites emerge from a family of rhombii solutions via a pitchfork bifurcation at $m = \widetilde{m}$.

\item[{\bf (iv)}]  For $m < -2$, there are four concave solutions corresponding to one geometrically distinct family of 
kites.

\item[{\bf (v)}]   For all other values of~$m$, there are no solutions.
\end{itemize}

\label{Thm:kites}
\end{theorem}

\pf
The formulas for $\rho$ and $\sigma$ come from solving equations (\ref{eq:rhoM}) and~(\ref{eq:sigmaM}).  Note that for fixed $m$, 
there are two possible values of $\rho$ depending on which sign is chosen.   Once $\rho$ is determined, there are two choices of $\sigma$
that arise from reflecting the kite about the $x$-axis.  Furthermore, we may interchange the values of $y_3$ and $y_4$ without changing the
values of $\rho$ and $\sigma$.  Thus, for a particular choice of $\rho$, there are four distinct kite relative equilibria described in abbreviated
coordinates by $(y_3, y_4)$, $(-y_3, -y_4)$, $(y_4, y_3)$, and $(-y_4, -y_3)$.   Translating and scaling each configuration so that $c=0$ and
$I=1$ gives four critical points of $H|_{\cal M}$.  All four solutions have the same shape and are equivalent under a reflection or a relabeling
of vortices 3 and~4.

To insure real solutions, we must have $\rho \in \mathbb{R}, \sigma^2 > 0,$ and the discriminant $\sigma^2 - 4\rho > 0$.  
It is straight-forward analysis to determine when these three conditions are met in terms of~$m$.  
This can be made rigorous using root-counting methods such as Sturm's theorem~\cite{sturm}.
The type of configuration (concave or convex) is governed by the sign of~$\rho$.  For example, if $0 < m < 1$, we obtain two distinct positive
values for $\rho$, each of which yield positive values for $\sigma^2$ and the discriminant.  Consequently, there are eight solutions 
and two geometrically distinct concave kite configurations.  At $m = 1$, we find that $\rho = 1$ is a double root of equation~(\ref{eq:rhoM}) and the two
families of kites merge into one solution given by $y_3 = \sqrt{3}, y_4 = 1/\sqrt{3}$ and its symmetric cousins.  The degeneracy of this particular configuration
is well known~\cite{albouy-four,MeyerSchmidt,palmore} (also see the second remark after Lemma~7.4 in~\cite{HRS}).

The other cases follow in a similar fashion.   When $+$ is chosen in equation~(\ref{eq:rhoForm}), the roots of $y^2 - \sigma y + \rho$ are complex for
$m > 1$ or $m < 0$.  At $m = 0$, the discriminant vanishes and we find $y_3 = y_4 = \sqrt{3}$, corresponding to a collision between vortices 3 and~4.  
On the other hand, taking $-$ in equation~(\ref{eq:rhoForm}) leads to several families of solutions.  Here, the solution at $m=0$ is $y_3 = 0, y_4 = \sqrt{3}$ and 
vortices 1, 2, and 3 are collinear while vortices 1, 2, and 4 form an equilateral triangle.  As $m$ decreases below 0, $\rho$ flips sign and the configuration
becomes convex.  The bifurcation at $\widetilde{m}$ comes from solving $\sigma^2 = 0$.  For this parameter value, we have $y_3 + y_4 = \sigma = 0$, so the configuration
is a rhombus.   As explained in~\cite{HRS}, the value $m = \widetilde{m}$ corresponds to a pitchfork bifurcation where the rhombus (two critical points
of $H|_{\cal M}$) bifurcates into the convex kites (four critical points of $H|_{\cal M}$).  This follows by replacing $m$ with $1/m$ and relabeling
the vortices so that our kite configurations match the framework used in~\cite{HRS}.
\enpf

\begin{remark}
We note that the results of Theorem~\ref{Thm:kites} agree with Table~1 in~\cite{HRS} once the transformation $m \mapsto 1/m$ is applied to convert
a $\mbox{Kite}_{12}$ configuration into a $\mbox{Kite}_{34}$. 
\end{remark}

\subsection{Two pairs of equal-strength vortices: asymmetric family}
\label{SubSec:Asymm}

For any $m \in (-1, 1)$, there exists a one-parameter family of asymmetric four-vortex relative equilibria with 
circulations $\Gamma_1 = \Gamma_2 = 1$ and $\Gamma_3 = \Gamma_4 = m$.   The existence of this family
was proven in~\cite{HRS}.  Because of the asymmetry, explicit
formulas for the positions of the vortices are more difficult to establish than with the kite families of the
preceding section.   We use Cartesian coordinates and the well-known Dziobek equations~\cite{dzio}, from which it is possible
to compute a Gr\"{o}bner basis.

\begin{figure}[tb]
\begin{center}
\includegraphics[width=268bp]{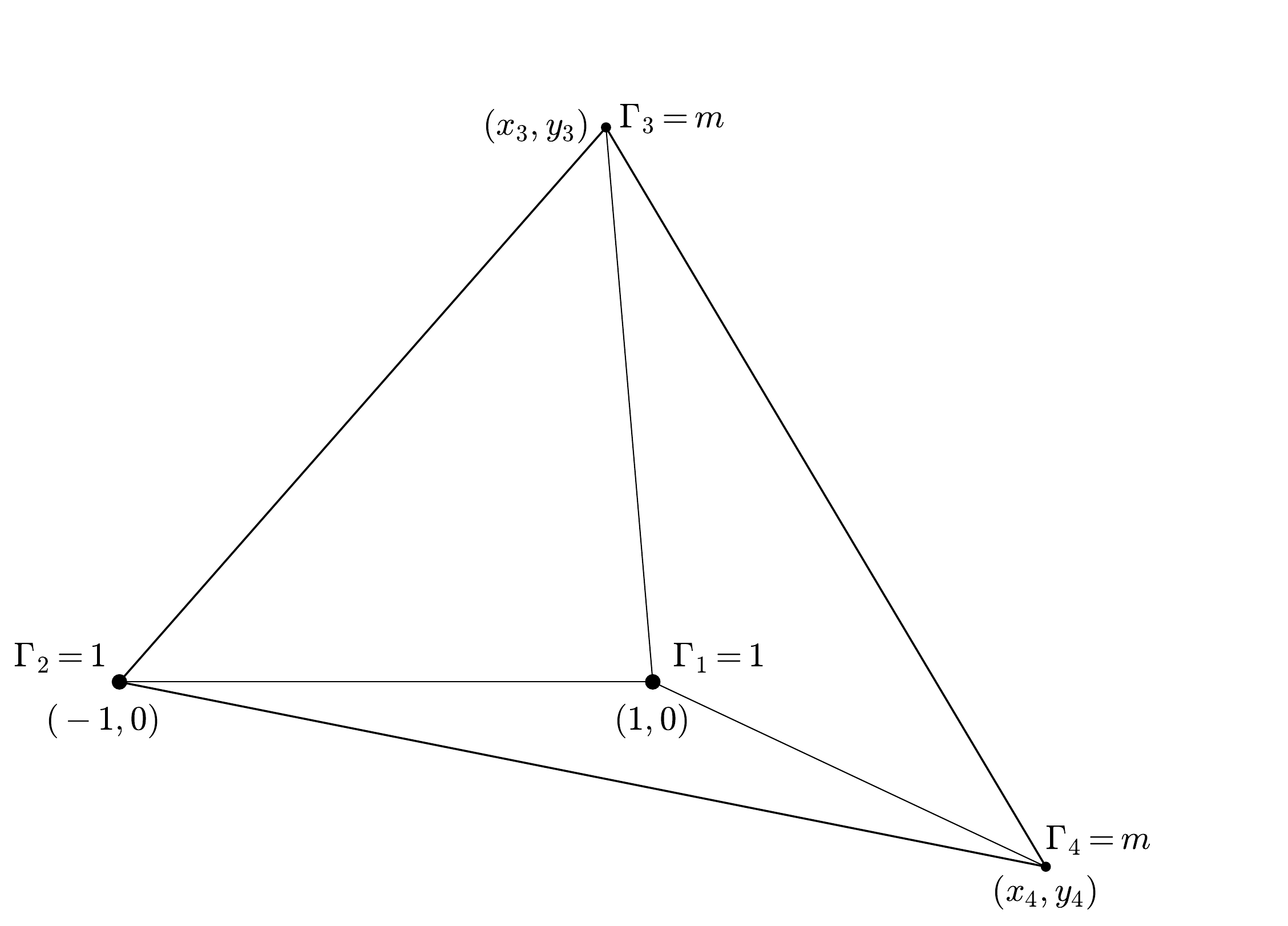}
\hspace{-0.35in}
\includegraphics[width=268bp]{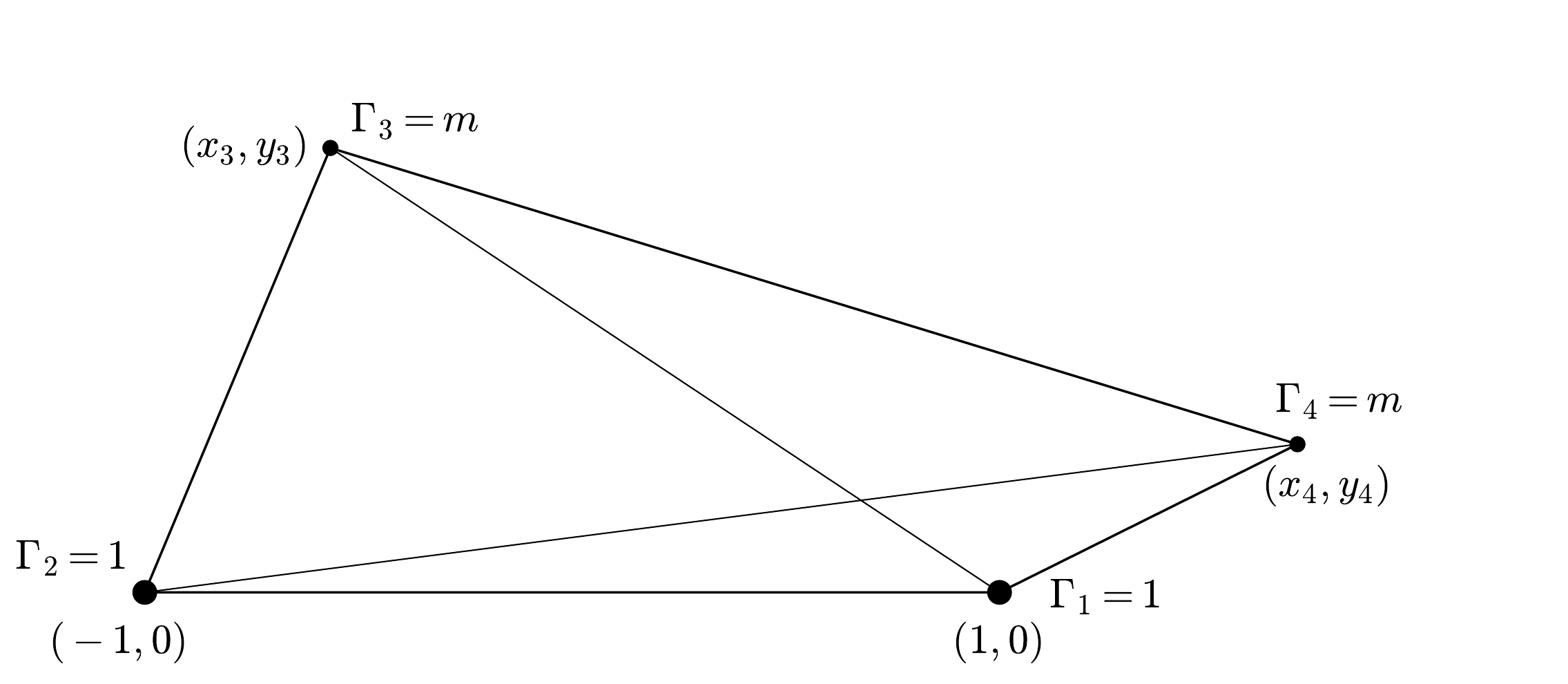}
\end{center}
\vspace{-0.2in}
\caption{Concave and convex asymmetric relative equilibria for the parameter values $m = 0.6$ (left) and $m = -0.6$ (right).}
\label{fig:Asymm}
\end{figure}

Let the positions of the vortices be $z_1 = (1, 0), z_2 = (-1, 0), z_3 = (x_3, y_3)$, and $z_4 = (x_4, y_4)$, where $x_3, x_4, y_3, y_4$
are unknown (see Figure~\ref{fig:Asymm}).  The Dziobek equations are derived using the mutual distances~$r_{ij}$ as variables, where an extra condition
(the vanishing of the Cayley-Menger determinant) is required to insure the vortices lie in the plane~\cite{sch, HRS}.
Introduce the variables $s_{ij} = r_{ij}^2$ and $\tau_{ij} = r_{ij}^{-2}$, and let $\Delta_i$ represent the oriented area
of the triangle formed by all of the vortices except for the $i$th vortex.   Assuming the $r_{ij}$ variables describe an actual configuration in
the plane, the following four equations are necessary and sufficient for a four-vortex relative equilibrium:
\begin{eqnarray}
\frac{\Gamma_1 \Delta_2}{\Gamma_2 \Delta_1}  & = &  \frac{\tau_{23} - \tau_{24}}{\tau_{13} - \tau_{14}} , \label{eq:Dzio1} \\[0.07in]
\frac{\Gamma_1 \Delta_3}{\Gamma_3 \Delta_1}  & = &  \frac{\tau_{23} - \tau_{34}}{\tau_{12} - \tau_{14}} ,  \label{eq:Dzio2} \\[0.07in]
\frac{\Gamma_3 \Delta_4}{\Gamma_4 \Delta_3}  & = &  \frac{\tau_{14} - \tau_{24}}{\tau_{13} - \tau_{23}} , \label{eq:Dzio3} \\[0.1in]
(s_{13} - s_{12})(s_{23} - s_{34})(s_{24} - s_{14}) & = &  (s_{12} - s_{14})(s_{24} - s_{34})(s_{13} - s_{23}).  \label{eq:Dzio4}
\end{eqnarray}
The signed areas $\Delta_i$ satisfy $\sum_i \Delta_i = 0$.   In our coordinates they are given by
\begin{eqnarray*}
\Delta_1 \; = \;   \frac{1}{2}\left( x_4 y_3 - x_3 y_4 + y_3 - y_4 \right),  & \quad &  \Delta_3 \; = \;  y_4,  \\[0.07in]
\Delta_2 \; = \;   \frac{1}{2}\left( x_3 y_4 - x_4 y_3 + y_3 - y_4 \right),  & \quad &  \Delta_4 \; = \;  -y_3. 
\end{eqnarray*}
Substituting the expressions for $\Delta_i$ into equations (\ref{eq:Dzio1}), (\ref{eq:Dzio2}), and~(\ref{eq:Dzio3}) yields 
\begin{eqnarray}
y_3 s_{14} s_{24} (s_{23} - s_{13}) + y_4 s_{13} s_{23} (s_{24} - s_{14}) & = & 0  ,  \label{eq:Asymm1} \\
s_{23} s_{24} (s_{14} - s_{13})(x_3 y_4 - x_4 y_3 + y_3 - y_4) + s_{13} s_{14} (s_{23} - s_{24})(x_4 y_3 - x_3 y_4 + y_3 - y_4) & = & 0 ,  \label{eq:Asymm2}  \\ 
2m s_{14} (s_{34} - s_{23})(x_4 y_3 - x_3 y_4 + y_3 - y_4) - y_4 s_{23} s_{34} (s_{14} - 4) & = &  0.  \label{eq:Asymm3}
\end{eqnarray}

Using $x_3, x_4, y_3, y_4$ and the unknown $s_{ij}$ as variables, it is possible to compute a Gr\"{o}bner basis for the polynomial ideal ${\cal A_S}$ generated by equations 
(\ref{eq:Dzio4}), (\ref{eq:Asymm1}),~(\ref{eq:Asymm2}) and the five distance relations involving the~$s_{ij}$ (e.g., $s_{13} = (x_3 - 1)^2 + y_3^2$).
To exclude the symmetric solutions, we saturate with respect to $s_{13} - s_{24}$,
$s_{14} - s_{23}$, $s_{13} - s_{14}$, $s_{23} - s_{24}$, $s_{13} - s_{23}$, $s_{14} - s_{24}$, $y_3$, and $y_4$.  
The computation was performed using the software Magma and took 5.6 minutes of CPU time.
With a lex ordering that eliminates all variables except $y_3$ and $y_4$, we obtain the polynomial $N_1 \cdot N_2$, where
$$
N_1 \; = \;  4y_3^2 y_4^2 - 3(y_3 - y_4)^2 
\quad \mbox{ and} \quad
N_2 \; = \;  (y_3^2 + y_4^2)^2 - 3(y_3 - y_4)^2.
$$

For $-1 < m < 1$, there are no real solutions that satisfy both $N_1 = 0$ and equations (\ref{eq:Dzio4}) through~(\ref{eq:Asymm3}).   
To see this, we append equation~(\ref{eq:Asymm3}) and $N_1 = 0$ to
the ideal ${\cal A_S}$ (along with the previous saturations) 
and compute a lex Gr\"{o}bner basis with an ordering of the form $\cdots > y_3 > s_{34} > m > y_4$.  This leads to the simple expression 
$s_{34} = 4m/(2m + 1)$.   Next, adding this new condition along with the previous two to ${\cal A_S}$ and 
computing a lex Gr\"{o}bner basis with an ordering of the form $\cdots > y_3 > y_4 > s_{34} > m > s_{14}$, we obtain the polynomial
$$
N_3 = (m + \casefrac{1}{2} )^2 x^4 - (10m^2 + 11m + 3)x^3 + (32m^2 + 42m + 14)x^2 - (40m^2 + 64m + 24)x +16(m+1)^2 ,
$$
where $x = s_{14}$.  Replacing $x$ by $2/x$ in $N_3$ and multiplying through by $\casefrac{1}{4} x^4$ yields
the polynomial~$P_1$ discussed in Section~6.1 of~\cite{HRS}.  Due to Lemma~6.1 in~\cite{HRS}, for any choice of
$m \in (-1, 1)$, $N_3$ has no positive roots in~$x$.  Consequently, there are no relative equilibria solutions corresponding
to solutions of $N_1 = 0$.

The asymmetric relative equilibria for $-1 < m < 1$ arise from solutions to $N_2 = 0$.  After adding equation~(\ref{eq:Asymm3}) and $N_2 = 0$ to
the ideal ${\cal A_S}$, we compute a lex Gr\"{o}bner basis with an ordering of the form $\cdots > y_3 > y_4 > m > s_{34}$.     
Using Magma, this computation took 19.4 minutes of CPU time and produced a basis with 66 elements.  The last polynomial
in this basis gives $s_{34} = 4(m + 2)$.   Another useful polynomial in the basis is $y_3 y_4 + \casefrac{3}{32} s_{34}^2 - \casefrac{9}{8} s_{34} + 3$.
By appending these two conditions to our ideal, we quickly compute two
lex Gr\"{o}bner basis (under 30 seconds of CPU time each):
one that eliminates all variables except for $s_{14}, x_3, x_4, $ and~$m$, and another that
eliminates all variables except for $y_3, y_4, $ and~$m$.  These calculations yield
\begin{eqnarray}
x_3^2 + x_4^2 \; = \; 3m + 5,  & \qquad  & y_3^2 + y_4^2 \; = \; 3(m + 1),  \label{eq:x3xrSqr} \\[0.05in]
2x_3 x_4 \; = \; m(3m + 5),   &  \qquad  &  2y_3 y_4 \;  = \;  -3m(m + 1),  \label{eq:x3x4Prod}
\end{eqnarray}
and 
\begin{equation}
(3m + 5)(s_{14} + 2x_4 - 3m - 4) - 2x_4^2 \; = \; 0,
\label{eq:s14Cond}
\end{equation}
from which we can derive formulas for the unknown position variables.

\begin{theorem}
Consider the four-vortex problem with vorticities $\Gamma_1 = \Gamma_2 = 1$ and $\Gamma_3 = \Gamma_4 = m$, with $-1 < m < 1$.
The vector $z = (1,0, -1,0, x_3, y_3, x_4, y_4)$ gives an asymmetric relative equilibrium, where 
\begin{eqnarray}
x_3 \; = \;  \frac{1}{2} \sqrt{3m + 5} \left( \sqrt{m+1} - \sqrt{1 - m} \, \right),   & \qquad  &  y_3 \; = \;  \frac{\sqrt{3}}{2} \left( m + 1 + \sqrt{1 - m^2} \, \right),  \label{eq:Asymmx3y3} \\[0.07in]
x_4 \; = \;  \frac{1}{2} \sqrt{3m + 5} \left( \sqrt{m+1} + \sqrt{1 - m} \, \right),  & \qquad  &  y_4 \; = \;  -\frac{\sqrt{3}}{2} \left( m + 1  -  \sqrt{1 - m^2} \, \right).   \label{eq:Asymmx4y4} 
\end{eqnarray}
For $0 < m < 1$, the configuration is concave with vortex~1 in the interior, while
for $-1 < m < 0$, the configuration is convex with adjacent equal pairs of vortices.   
The transition from concave to convex configuration at $m = 0$ is smooth.  
All solutions, regardless of the value of~$m$, have the same angular velocity $\omega = 1/2$.  
For each $m \in (-1, 1)$, there are a total of eight asymmetric relative equilibria, each geometrically equivalent to the given solution~$z$.  
\label{Thm:AS}
\end{theorem}

\pf
Equations (\ref{eq:x3xrSqr}) and~(\ref{eq:x3x4Prod}) imply 
\begin{eqnarray*}
(x_3 - x_4)^2 \; = \; (3m + 5)(1 - m),  & \qquad  & (y_3 - y_4)^2  \; = \; 3(m + 1)^2,   \\[0.05in]
(x_3 + x_4)^2 \; = \; (3m + 5)(m + 1),   &  \qquad  &  (y_3 + y_4)^2 \;  = \;  3(1 - m^2),
\end{eqnarray*}
which, in turn, lead to the expressions for $x_3, x_4, y_3$, and $y_4$ stated in the theorem.    The signs in front of $\sqrt{1-m}$ and $\sqrt{1-m^2}$ must be chosen carefully
to insure that equation~(\ref{eq:s14Cond}) is also satisfied.   Using Maple, equations (\ref{eq:Asymmx3y3}) and~(\ref{eq:Asymmx4y4}) were confirmed to satisfy 
system~(\ref{eq:rel-equ}), the equations for a relative equilibrium in Cartesian coordinates.

For $0 < m < 1$, we have $x_3, y_3, x_4 > 0$ and $y_4 < 0$, so vortices 3 and~4 reside in quadrants I and~IV, respectively.  To prove that
the configuration is concave, we show that the $x$-intercept of the line between vortices 3 and~4 is greater than~1.
This implies that vortex~1 lies inside the triangle formed by the other three vortices.  The $x$-intercept of the line between $z_3$ and $z_4$ is
$(x_4 y_3 - x_3 y_4)/(y_3 - y_4)$, so it suffice to show that
\begin{equation}
x_4 y_3 - x_3 y_4 + y_4 - y_3 \; > \; 0.
\label{eq:ConcvPf}
\end{equation}
Substituting in equations (\ref{eq:Asymmx3y3}) and~(\ref{eq:Asymmx4y4}), inequality~(\ref{eq:ConcvPf}) reduces to
$$
\sqrt{3(m + 1)} \left(  \sqrt{3m + 5} - \sqrt{m + 1} \, \right) \; > \; 0,
$$
which is clearly valid for $0 < m < 1$.

At $m = 0$, $z_3 = (0, \sqrt{3})$ and $z_4 = (\sqrt{5}, 0)$, so vortices 2, 1, and~4 are collinear and
vortices 2, 1, and~3 form an equilateral triangle.  As $m$ becomes negative, both $x_3$ and $y_4$ flip signs, moving vortices 3 and~4 into
quadrants II and~I, respectively.  The configuration is now convex because inequality~(\ref{eq:ConcvPf}) still holds, so
the $x$-intercept of the line between vortices 3 and~4 remains larger than~1.  
The transition from concave to convex configuration is smooth
because the derivatives of $x_3, x_4, y_3$, and $y_4$ with respect to~$m$ evaluated at $m = 0$ exist and are all nonzero.

Formulas for the $s_{ij}$ follow easily from equations (\ref{eq:Asymmx3y3}) and~(\ref{eq:Asymmx4y4}).  By straight-forward analysis, we find that 
when $0 < m < 1$ (concave), $s_{24} > s_{34} > s_{23} > s_{13} > s_{12}  > s_{14}$, while for the case $-1 < m < 0$ (convex), we have
$s_{24} > s_{34} > s_{12} > s_{13} > s_{23} > s_{14}.$   These inequalities are all strict except for $m = 0$ or $m = 1$, which serves to verify
the asymmetry of the configuration.  Note that at $m = 1$, we find $z_3 = (2, \sqrt{3})$ and $z_4 = (2, -\sqrt{3})$, so the outer triangle
is equilateral with vortex~1 at the center.  Thus, as with the kite solutions of Theorem~\ref{Thm:kites},
the asymmetric family bifurcates out of the degenerate equilateral triangle solution.

The angular impulse $I$ (with respect to the center of vorticity) can be written nicely in terms of the $s_{ij}$ as
\begin{equation}
I \; = \;  \frac{1}{\Gamma} \sum_{i < j} \Gamma_i  \Gamma_j  r_{ij}^2  \; = \;   \frac{1}{\Gamma} \sum_{i < j} \Gamma_i  \Gamma_j  s_{ij} \, .
\label{eq:InertiaMD}
\end{equation}
Using (\ref{eq:Asymmx3y3}) and~(\ref{eq:Asymmx4y4}) implies that $I = 2(m^2 + 4m + 1) = 2L$, so that $\omega = L/I = 1/2$ for any value of~$m$.

For each $m \in (-1, 1)$, there are four distinct asymmetric relative equilibria given in abbreviated
coordinates by $(x_3, y_3, x_4, y_4)$, $(x_3, -y_3, x_4, -y_4)$, $(x_4, y_4, x_3, y_3)$, and $(x_4, -y_4, x_3, -y_3)$.   
An additional four solutions are generated by interchanging $z_1$ and $z_2$ in each of these configurations.
Thus, after translating and scaling each configuration so that $c=0$ and
$I=1$,  we obtain eight critical points of~$H|_{\cal M}$.  
All eight solutions have the same shape and are equivalent under a reflection or a relabeling
of vortices 3 and~4, or 1 and~2, or both.  

\enpf

\subsection{Applying the Morse inequalities}
\label{SubSec:ApplyingMorse}

We now apply the Morse inequalities to the planar four-vortex problem with two pairs of equal vorticities.  The first step is
to verify that $H|_{\cal M}$ is a Morse function by checking that all critical points are nondegenerate.  This can be accomplished using
Theorems \ref{Thm:kites} and~\ref{Thm:AS} and Gr\"{o}bner bases.
Instead of translating and rescaling the kite and asymmetric solutions so that $c = 0$ and $I = 1$, it is easier to work with the coordinates and formulas
given in Theorems \ref{Thm:kites} and~\ref{Thm:AS}.  If $z$ is nondegenerate in this setting, then the corresponding configuration in~${\cal N}$
is also nondegenerate.

For a relative equilibrium $z$ with angular velocity~$\omega$, the modified Hessian is given by the matrix $M^{-1} D^2H(z) + \omega I$.  
Regardless of the signs of the circulations $\Gamma_i$, a vector in the kernel of the Hessian $D^2G(z)$ will also be in the kernel
of the modified Hessian.  Following the
arguments of Section~\ref{subsec:trivialEvals}, the modified Hessian always has the trivial eigenvalues $\omega, \omega, 2\omega,$ and $0$, with corresponding
eigenvectors $s, Ks, z$ and $Kz$, respectively.   To verify that $z$ is nondegenerate, we must show that all of the other eigenvalues are nonzero.   
This is equivalent to showing the nontrivial eigenvalues of $M^{-1} D^2H(z)$ are not equal to~$-\omega$.  The following lemma, which also applies
in the case of mixed-sign vorticities, gives a necessary and sufficient condition for nondegeneracy in terms of the coefficients of the characteristic
polynomial of $M^{-1} D^2H(z)$.

\begin{lemma}
Suppose that $z$ is a relative equilibrium of the four-vortex problem and let 
$$
R(\mu) \; = \;  \mu^8 + c_6 \mu^6 + c_4 \mu^4 + c_2 \mu^2
$$
be the characteristic polynomial of the matrix $M^{-1} D^2H(z)$.  Then, $z$ is nondegenerate if and only if
\begin{equation}
c_4 + 2 \omega^2 c_6 + 3 \omega^4 \;  \neq \;  0  .
\label{eq:nondeg}
\end{equation}
\end{lemma}

\pf
Recall that identity~(\ref{eq:anti-commute}) implies that if $v$ is an eigenvector of $M^{-1} D^2H(z)$ with eigenvalue~$\mu$, then
$Kv$ is an eigenvector with eigenvalue $-\mu$.  Thus, eigenvalues of $M^{-1} D^2H(z)$ come in pairs of the form~$\pm \mu_j$.
For mixed-sign circulations, the $\mu_j$ may be complex.  The trivial eigenvalues are $0, 0, \pm \omega$.  
Denote the remaining four eigenvalues as $\pm \mu_1, \pm \mu_2$.  If the eigenvalues are complex, we have $\mu_2 = \overline{\mu_1}$.
Expanding $R(\mu)$, we find that
\begin{eqnarray*}
R(\mu)  & = &  \mu^2 (\mu^2 - \omega^2) (\mu^2 - \mu_1^2) (\mu^2 - \mu_2^2)  \\[0.05in]
&  =  &  \mu^8  - (\mu_1^2 + \mu_2^2 + \omega^2) \mu^6  + [\mu_1^2 \mu_2^2 + \omega^2(\mu_1^2 + \mu_2^2)]  \mu^4 - \omega^2 \mu_1^2 \mu_2^2 \, \mu^2 .
\end{eqnarray*}
Therefore, $c_6 = -(\mu_1^2 + \mu_2^2 + \omega^2)$ and $c_4 = \mu_1^2 \mu_2^2 + \omega^2(\mu_1^2 + \mu_2^2)$, which in turn, yields
$$
c_4 + 2\omega^2 c_6 + 3 \omega^4 \; = \;  (\mu_1^2 - \omega^2) (\mu_2^2 - \omega^2).
$$
Since $z$ is nondegenerate if and only if $\pm \mu_i \neq -\omega$ for each~$i$, the result follows.
\enpf

\begin{lemma}
For the planar four-vortex problem with vortex strengths $\Gamma_1 = \Gamma_2 = 1$ and $\Gamma_3 = \Gamma_4 = m$, 
the Hamiltonian $H$ restricted to the quotient manifold~${\cal M}$ is a Morse function for each $m$ satisfying $0 < m < 1$.
\label{lemma:Morse}
\end{lemma}

\pf
Of the 34 relative equilibria, the 12 collinear configurations are known to be nondegenerate~\cite{palmore, rick-book}
when $m > 0$.  The 6 convex configurations (rhombii and isosceles trapezoids) were shown to be nondegenerate for $0 < m \leq 1$ in~\cite{g:stability}.
It remains to check the nondegeneracy of the 16 concave kite and asymmetric configurations discussed in Sections \ref{SubSec:Kites}
and~\ref{SubSec:Asymm}.

We will verify condition~(\ref{eq:nondeg}) holds for any choice of $m \in (0,1)$.  The values of the coefficients $c_4$ and $c_6$ can be
expressed in terms of the traces of the powers of $M^{-1} D^2H(z)$ using the Leverrier-Souriau-Frame algorithm (see p. 504 in~\cite{CarlMeyer}).
Set $C = M^{-1} D^2H(z)$ and let  tr$(\ast)$ denote the trace of a matrix.  We find that
\begin{eqnarray*}
c_6 & = &  \frac{1}{2}  \left(  [\mbox{tr} (C) ]^2 - \mbox{tr} (C^2)  \right), \\[0.07in] 
c_4 & = &  \frac{1}{24}\left[ \mbox{tr} (C) \right]^4 + \frac{1}{3} \mbox{tr} (C) \cdot \mbox{tr} (C^3)  - \frac{1}{4} [\mbox{tr} (C)]^2 \cdot \mbox{tr} (C^2)
+ \frac{1}{8} [\mbox{tr} (C^2) ]^2 - \frac{1}{4} \mbox{tr} (C^4).
\end{eqnarray*}
Although these formulas appear daunting, they are greatly simplified in our setting because the characteristic polynomial of $M^{-1} D^2H(z)$
is even.  Consequently, tr$(C)=0$ and we have
$$
c_6 \; = \;  -\frac{1}{2} \mbox{tr}(C^2)  \quad  \mbox{and} \quad
c_4 \; = \;  \frac{1}{8} [\mbox{tr} (C^2) ]^2 - \frac{1}{4} \mbox{tr} (C^4).
$$
Formulas for $c_4$ and $c_6$ in terms of the positions $(x_i, y_i)$, mutual distances $r_{ij}$, and vorticities~$\Gamma_i$ are lengthy and
included in the appendix.

We begin with the kite solutions of Theorem~\ref{Thm:kites}.  Substituting the formulas for the positions directly into~(\ref{eq:nondeg}) is too cumbersome;
there are 5,289 terms in the numerator.  Instead, we use Gr\"{o}bner bases to eliminate all variables except for~$m$.   Let ${\cal P}_K$ be the polynomial
ideal generated by the numerator of $c_4 + 2 \omega^2 c_6 + 3 \omega^4$ evaluated at $z = (1,0, -1,0, 0, y_3, 0, y_4)$, along with
equations (\ref{eq:sigma}) through~(\ref{eq:sigmaM}).   The value of $\omega$ is found from $\omega = L/I$, where $I$ is given by equation~(\ref{eq:InertiaMD}).
We saturate with respect to $y_3 + y_4$ to eliminate the rhombus solutions.  
The variety of ${\cal P}_K$ will contain all of the degenerate kite configurations.
Computing a Gr\"{o}bner basis for ${\cal P}_K$ with respect to the lex order 
$y_3 > y_4 > \sigma > \rho > m$ yields the polynomial
\begin{equation}
m^4 (m - 1)(m^2 + 4m +1)^5(2m + 1)(m^2 - 4m - 9)(3m + 5)(m + 2)^{10}.
\label{eq:kitesMpoly}
\end{equation}
Note that $m = 1$ and $m = -2 \pm \sqrt{3}$ (where $L = 0$), which are known degenerate cases, are roots of this polynomial.  Since~(\ref{eq:kitesMpoly})
has no roots strictly between 0 and~1, condition~(\ref{eq:nondeg}) is satisfied for any kite configuration with $m \in (0, 1)$.

For the asymmetric family, we have $\omega = 1/2$ for any~$m$.   Introduce the auxiliary variables 
$u_1 = \sqrt{3m + 5} \, , u_2 = \sqrt{m + 1} \,$, and $u_3 = \sqrt{1-m} \,$, and let $p_{dg}$ be the numerator 
of $16c_4 + 8c_6 + 3$ evaluated at $z = (1,0, -1,0, x_3, y_3, x_4, y_4)$ using
formulas (\ref{eq:Asymmx3y3}) and~(\ref{eq:Asymmx4y4}) with the $u_i$ variables.  
Define ${\cal P}_A$ to be the polynomial ideal generated by $p_{dg}$ and the three relations $u_1^2 - (3m+5),
u_2^2 - (m+1)$, and $u_3^2 - (1 - m)$.  
The variety of ${\cal P}_A$ contains the degenerate solutions from the asymmetric family of relative equilibria.  
Computing a Gr\"{o}bner basis for ${\cal P}_A$ with respect to a lex order that eliminates all variables except for~$m$ 
produces the polynomial
\begin{equation}
(m + 2)^9 (3m + 5) (m + 1)^8 (m^2 + 4m + 1) (m - 1) .
\label{eq:AsymmMpoly}
\end{equation}
As expected, we recover the degeneracy of the equilateral triangle solution at $m = 1$, as well as the case when $L = 0$.
Since~(\ref{eq:AsymmMpoly}) has no roots for $0 < m < 1$, the asymmetric family is nondegenerate for these parameter values and the
proof is complete.
\enpf

\begin{remark}
\begin{enumerate}
\item  From the proof of Lemma~\ref{lemma:Morse}, we see that the kite and asymmetric solutions, when they exist, are also nondegenerate for $m < 0$,
except when $m^2 + 4m + 1 = 0$.  However, we cannot conclude that $H|_{\cal M}$ is a Morse function for $m < 0$ because
we do not know that the collinear relative equilibria are also nondegenerate.  

\item  Remarkably, even though the polynomials arising from~(\ref{eq:nondeg}) contain thousands of terms ($p_{dg}$ has $11,\!644$ terms), 
the Gr\"{o}bner basis calculations producing polynomials (\ref{eq:kitesMpoly}) and~(\ref{eq:AsymmMpoly}) each take under a second of CPU time using Maple.
\end{enumerate}
\end{remark}

\begin{theorem}
Consider the four-vortex problem with vortex strengths $\Gamma_1 = \Gamma_2 = 1$ and $\Gamma_3 = \Gamma_4 = m$, where $m \in (0,1)$ is
a parameter.   The concave relative equilibria (the kite and asymmetric families) are unstable for all values of~$m$.  For each relative equilibrium, the four
nontrivial eigenvalues consist of one pair of real values $\pm \lambda$ and one pair of pure imaginary values $\pm i \beta$, with $|\beta| \leq \omega$.
\end{theorem}

\pf
From equation~(\ref{eq:Poincare}), the Poincar\'{e} polynomial for the case $n=4$ is $P(t) = 1 + 5t + 6t^2$.  Thus, the Morse inequalities can be written as
\begin{equation}
\gamma_0 + \gamma_1 t + \gamma_2 t^2 \; = \;  1 + 5t + 6t^2 + (1 + t)(r_0 + r_1 t),
\label{eq:Morse}
\end{equation}
where $r_0$ and $r_1$ are non-negative integers.

Based on the work in~\cite{HRS}, there are exactly 34 relative equilibria for each value of $m$.  The six convex configurations are shown to be linearly
stable (and therefore minima) for $m > 0$ in~\cite{g:stability}, so $\gamma_0 \geq 6$.  Since the 12 collinear solutions have index 2, we also know
that $\gamma_2 \geq 12$.   Equating the coefficients of the constant and quadratic terms on each side of equation~(\ref{eq:Morse}) gives $r_0 \geq 5$ and $r_1 \geq 6$. 
On the other hand,  setting $t = 1$ in~(\ref{eq:Morse}) gives
$$
34 \; = \;  \gamma_0 + \gamma_1 + \gamma_2 \;  = \;  12 + 2(r_0 + r_1) \quad \Longrightarrow  \quad  r_0 + r_1 = 11 .
$$
It follows that $r_0 = 5$ and $r_1 = 6$.  This, in turn, yields $\gamma_0 = 6, \gamma_1 = 16$, and $\gamma_2 = 12$.  
Thus, the 16 remaining relative equilibria (the concave solutions) each have index one and by Theorem~\ref{Thm:Main}, one pair of real eigenvalues and one pair
of pure imaginary eigenvalues.  The estimate for $\beta$ comes from the third point in Remark~\ref{remark:mainThm}.
\enpf

\begin{remark}
The eigenvalues for the kite and asymmetric families were computed numerically for specific cases when $m < 0$.  
None of the solutions were linearly stable.  For the kite configurations, the nontrivial eigenvalues contained either 
one or two real pairs, while the asymmetric family had either one real pair or a complex quartuplet.  Changes in the eigenvalue structure occurred at
$m = -2 \pm \sqrt{3}$, the two values for which $L = 0$.  There did not appear to be any connection between the Morse index and the number of real
or complex eigenvalues when $m < 0$.  
\end{remark}

\section{Conclusion}

We have taken a Morse theoretical approach to study the stability of relative equilibria in the
planar $n$-vortex problem with positive circulations.  
Treating relative equilibria as critical points of a smooth function on a manifold, we have shown that the Morse index 
is equivalent to the number of pairs of real (nonzero) eigenvalues.  
In essence, the greater the index, the more unstable the relative equilibrium becomes.

For a fixed choice of positive circulations, we have shown that relative equilibria cannot
accumulate on the collision set.   This allows us to restrict the space to a compact set and ensures,
assuming nondegeneracy, that the number of critical points is finite.  The Morse inequalities
were utilized in the two equal pairs problem to show that the concave kite and asymmetric
families of relative equilibria are unstable, each with one real pair of eigenvalues.
The most difficult part of the calculation was verifying that these solutions were nondegenerate.

In future work, we hope to apply these same techniques to relative equilibria of the four-vortex problem with three
equal circulations (e.g., $\Gamma_1 = \Gamma_2 = \Gamma_3 = 1$, $\Gamma_4 = m$).
The collinear solutions for this case have been rigorously explored in~\cite{gBrian}, with 
linearly stable solutions located for $m$ close to~$-1$.  It would also be interesting to 
extend the theory developed in this work to the case of mixed-sign circulations.
There, the level surface $I = 1$ becomes a hyperboloid and the topology changes dramatically.  
The circulation matrix $M$ is no longer positive definite, so the key matrix $M^{-1} D^2H(x)$ 
may have complex eigenvalues and the nice factorization of the characteristic polynomial in~(\ref{eq:charPoly})
is lost.  Nevertheless, it may be possible to apply other index theories (e.g., the Maslov index) in this setting
to obtain fruitful results~\cite{long, BJP2}.

\section{Appendix}

Here we provide formulas for the key coefficients $c_4$ and $c_6$, obtained from 
$$
c_4 \; = \;  \frac{1}{8} [\mbox{tr} (C^2) ]^2 - \frac{1}{4} \mbox{tr} (C^4)  \quad  \mbox{and} \quad
c_6 \; = \;  -\frac{1}{2} \mbox{tr}(C^2) ,
$$
where $C = M^{-1} D^2H(z)$.  For the case $n=4$, $C$ is an $8 \times 8$ matrix of the form
$$  
\begin{bmatrix}
C_{11}  &  C_{12}  &  C_{13}  &  C_{14} \\
C_{21}  &  C_{22}  &  C_{23}  &  C_{24} \\
C_{31}  &  C_{32}  &  C_{33}  &  C_{34} \\
C_{41}  &  C_{42}  &  C_{43}  &  C_{44} 
\end{bmatrix},
\quad \mbox{where }  \;
C_{ij}   =   \frac{1}{\Gamma_i} A_{ij} = 
\frac{\Gamma_j}{r_{ij}^4} 
\begin{bmatrix}
a_{ij} & b_{ij} \\[0.1in]
b_{ij} & -a_{ij}
\end{bmatrix}, \; 
C_{ii} =  - \sum_{j \neq i} C_{ij} \, ,
$$
and $a_{ij} = (y_i - y_j)^2 - (x_i - x_j)^2, b_{ij} = -2(x_i - x_j)(y_i - y_j)$.
Introduce the variables $s_{ij} = r_{ij}^2$, $r = \prod_{i < j} s_{ij}^2$,
$G_{ijk} = a_{ij} a_{ik} + b_{ij} b_{ik},$ and $H_{ijk} = a_{ij} b_{ik} - a_{ik} b_{ij}$.
The following identities are helpful in the calculation:
$$
C_{ji} \; = \;  \frac{\Gamma_i}{\Gamma_j} C_{ij},  \quad
C_{ij}^2 \; = \;  \frac{\Gamma_j^2}{s_{ij}^2} I_2,  \quad  \mbox{ and} \quad
C_{ij} C_{ji} \; = \;  \frac{\Gamma_i \Gamma_j}{s_{ij}^2} I_2 ,
\mbox{ where } 
I_2 = 
\begin{bmatrix}
1  &   0   \\
0  &   1
\end{bmatrix}.
$$
We find that
\begin{equation}
\mbox{tr}(C^2) \; = \;  2 \sum_{i < j} \frac{ (\Gamma_i + \Gamma_j)^2}{s_{ij}^2}
\, + \, 4 \sum_{i=1}^4  \sum_{\substack{ j < k \\[0.03in] j, k \neq i}}  \frac{ \Gamma_j \Gamma_k \, G_{ijk} }{s_{ij}^2 s_{ik}^2}  \, .
\label{eq:TraceSqrd}
\end{equation}

After a lengthy and tedious calculation involving the tr$(C^4)$ and the square of~(\ref{eq:TraceSqrd}), we compute (by hand) that
$
c_4  =   F_1 + 2F_2 + 2F_3 + 2F_4 + 4F_5 + 2F_6 + 4F_7,
$
where the $F_i$ are given as follows:
\begin{eqnarray*}
F_1  & = &  \sum_{i=1}^4  \sum_{\substack{ j < k \\[0.03in] j, k \neq i}}  \frac{ \Gamma_i^2 (\Gamma_i  + \Gamma_j + \Gamma_k)^2}{s_{ij}^2 s_{ik}^2}  +
 \sum_{\substack{ j = 2 \\[0.02in] k < l \\[0.03in] k \neq 1, j}}^4 \frac{ (\Gamma_1 + \Gamma_j)^2 (\Gamma_k + \Gamma_l)^2}{s_{1j}^2 s_{kl}^2} ,  \\[0.05in]
F_2  & = &  \sum_{ i < j}  \sum_{k \neq i, j}   \frac{ \Gamma_i \Gamma_j (\Gamma_i + \Gamma_j) (\Gamma_i + \Gamma_j + 2 \Gamma_k) G_{kij}}{s_{ij}^2 s_{ik}^2 s_{jk}^2} , \\
F_3  & = &  \sum_{i=1}^4  \sum_{\substack{ j < k \\[0.03in] j, k \neq i}}  \frac{ \Gamma_i^2 \Gamma_j \Gamma_k \, G_{ijk} }{s_{ij}^2 s_{ik}^2 s_{il}^2} , \\[0.05in]
F_4  & = &  \sum_{i < j}  \sum_{k \neq i, j}  \frac{ \Gamma_i \Gamma_k (\Gamma_i + \Gamma_l)^2 \, G_{jik}  + \Gamma_j \Gamma_l (\Gamma_j + \Gamma_k)^2 \, G_{ijl} }{s_{ij}^2 s_{il}^2 s_{jk}^2} , \\
F_5  & = &  \frac{1}{r} \sum_{i=1}^4  \sum_{\substack{ j < k \\[0.03in] j, k \neq i}}  s_{ij}^2 s_{ik}^2 \Gamma_i \Gamma_l \left[ \Gamma_k^2 \, G_{jlk} G_{lki} + \Gamma_j^2 \, G_{klj} G_{lji} 
+ \Gamma_j  \Gamma_k (G_{jlk}  G_{lki}  -  H_{jlk} H_{lki} ) \right] , \\
F_6  & = &  \sum_{i < j} \frac{\Gamma_i \Gamma_j}{s_{ij}^4} \left(  \, \sum_{k \neq i,j} \left[  \frac{\Gamma_k \, G_{ijk}}{s_{ik}^2} \left( \frac{ \Gamma_k \, G_{jik}}{s_{jk}^2} + \frac{ \Gamma_l \, G_{jil}}{s_{jl}^2} \right)
+ \frac{\Gamma_k \, H_{ijk}}{s_{ik}^2} \left( \frac{ \Gamma_k \, H_{jik}}{s_{jk}^2} + \frac{ \Gamma_l \, H_{jil}}{s_{jl}^2} \right) \right]  \right) , \\
F_7  & = &  \frac{1}{r} \sum_{\substack{ j = 2 \\[0.02in] k < l \\[0.03in] k \neq 1, j}}^4  s_{1j}^2 s_{kl}^2 \left[ \Gamma_1^2 \Gamma_j^2 \, G_{k1j} G_{l1j} + \Gamma_k^2 \Gamma_l^2 \, G_{1kl} G_{jkl} -
\Gamma_1 \Gamma_2 \Gamma_3 \Gamma_4 ( G_{k1j} G_{l1j} + H_{k1j} H_{l1j} )  \right]  .
\end{eqnarray*}
In any expression where indices $i, j, k, l$ appear together, it is always assumed that they are distinct.
For example, in $F_4$, if $i = 1, j = 4,$ and $k = 2$, then we take $l = 3$.  The number of terms (equivalent denominators) in each $F_i$ is  
$15, 4, 4, 12, 12, 24,$ and $3$, respectively, for a total of 74 terms in the coefficient~$c_4$.


\bibliographystyle{amsplain}

\end{document}